\documentclass[10pt,a4papers]{amsart}
\usepackage{amsmath}
\usepackage{amscd}
\usepackage{amsfonts}
\usepackage{amssymb}
\usepackage{ascmac}
\usepackage{latexsym}
\newtheorem{df}{Definition}[section]
\newtheorem{lm}[df]{Lemma}
\newtheorem{pr}[df]{Proposition}
\newtheorem{co}[df]{Corollary}
\newtheorem{theorem}[df]{Theorem}
\newtheorem{rk}[df]{Remark}
\newtheorem{ex}[df]{Example}
\newtheorem{prob}[df]{Problem}

\newtheorem{notation}[df]{Notation}

\newcommand{\bsquare}{\hbox{\rule{6pt}{6pt}}}
\begin{document}
\title{$D$-local algebraic fundamental groups of germs of normal complex
analytic singularities}\thanks{This paper grew out of the section 4 of the author's preprint \cite{ohno}.} 
\author{Koji Ohno
\\ Department of Mathematics Graduate School of Science \\Osaka University}
\date{\today}
\address[Koji Ohno]{Department of Mathematics Graduate School of Science, Osaka University, Toyonaka, Osaka, 560-0043, Japan}
\email[Koji Ohno]{koji@math.sci.osaka-u.ac.jp}
\maketitle 
\begin{abstract} In this paper, the notion of local algebraic fundamental groups of normal complex analytic singularities are generalized to certain profinite groups called $D$-local algebraic fundamental groups which turns out to be useful even for the study of usual local algebraic fundamental groups and the
Lefshetz type theorem for $D$-local algebraic fundamental groups is proved under certain conditions. The theorem yields, for example, the finiteness of the local algebraic
fundamental groups of a certain class of four dimensional singularities and will be useful for the classification of three dimensional purely log terminal singularities.
\end{abstract}
\tableofcontents
\begin{center}
Introduction
\end{center}
Reid-Wahl's cyclic covering trick such as taking canonical or log
canonical covers has been played the fundamental role in studying singularities. The underlying principle is that we can
understand singularities by taking suitable Galois covers. But in many cases, just like quotient singularity, or surface log
terminal singularity with branch loci (i.e., with a standard boudary in Shokurov's terminology ) as explained in
\cite{kns}, these cyclic covers is dominated by, in some sense, \lq\lq universal'' finite Galois
covers and singularities can be more easily seen if we can take this sort of coverings. In Section 1, we introduce the notion of $D$-{\it local algebraic fundamental groups}
(Definition~\ref{df:log pi1}) associated with a pair of a germ of complex singularity $X$ and a $\mbox{\boldmath $Q$}$-Cartier divisor $D$. In Section 2, We give a comparison theorem (Theorem~\ref{th:comparison}, 
which compare $D$-local fundamental groups with certain topological invariants which has been introduced in M. Kato
(\cite{kato}), M. Namba (\cite{namba}) and J.P. Serre (\cite{serre}, Appendix 6.4). In Section 3, we introduce {\it universal index one covers} (Definition~\ref{df:u-index
one cover}) associated with a pair of a germ of complex singularity $X$ and a $\mbox{\boldmath $Q$}$-Cartier divisor $D$ and deduce an exact sequence which relate the $D$-local fundamental group and the local fundamental group of its universal index one cover and obtain Lemma~\ref{lm:brieskorn's Lemma} which is a generalization of Brieskorn's lemma (\cite{brieskorn}, Lemma 2.6). In Section 4, we prove the Lefshetz type theorem for $D$-local
algebraic fundamental groups under certain conditions  (Theorem~\ref{th:pi1surj}) which enables us to study $D$-local
algebraic fundamental groups inductively on dimensions. 
\begin{center}
Notation and Terminology
\end{center}
Let $X$ be a normal Stein space or a germ of normal complex analytic spaces with a point $p\in X$. 
$\mbox{Weil }X$ is the free abelian group generated by prime divisors on $X$ and 
$\mbox{Div }X$ is the subgroup of $\mbox{Weil }X$ generated by Cartier divisors. 
An element of $\mbox{Weil }X$ (resp. $\mbox{Weil }X\otimes \mbox{\boldmath $Q$}$) 
is called an {\it integral divisor} 
(resp. a {\it {\bf Q}-divisor}). $\mbox{Div}_{\mbox{\boldmath $Q$}}X$ is the 
{\bf Q}-submodule of $\mbox{Weil }X\otimes \mbox{\boldmath $Q$}$ generated by $\mbox{Div }X$. 
An elements $D$ of $\mbox{Weil }X\otimes \mbox{\boldmath $Q$}$ is said to be {\it {\bf Q}-Cartier}, 
if $D\in \mbox{Div}_{\mbox{\boldmath $Q$}}X$. For $D\in \mbox{Div}_{\mbox{\boldmath $Q$}}X$, the index 
$[\mbox{\boldmath $Z$}D:\mbox{\boldmath $Z$}D\cap \mbox{Div }X]\in \mbox{\boldmath $N$}$ is called the index of $D$
at $p\in X$ denoted by $\mbox{ind}_pD$. Let $f:Y\rightarrow X$ be a finite morphism between 
normal Stein spaces or germs of normal complex analytic spaces. The pull-back homomorphism 
$f^{\ast}:\mbox{Div }X\rightarrow \mbox{Div }Y$ canonically  extends to a homomorphism 
$f^{\ast}:\mbox{Div}_{\mbox{\boldmath $Q$}}X\rightarrow \mbox{Div}_{\mbox{\boldmath $Q$}}Y$. We note that canonical divisor on 
the smooth loci of $X$ extends to a divisor $K_X$ on $X$ by the Remmert-Stein's theorem (see, for example, \cite{gunning-rossi}, Chap. V, Sec. D, Theorem 5).  Let
$\Delta$ be a {\bf Q}-divisor on $X$ such that $K_X+\Delta$ is {\bf Q}-Cartier.
$\Delta_Y$ is a {\bf Q}-divisor on $Y$ defined by $\Delta_Y:=f^{\ast}(K_X+\Delta)-K_Y$. {\bf Q}-divisor $\Delta$ 
is called {\it {\bf Q}-boundary}, if $\mbox{mult}_{\Gamma}\Delta\in [0,1]\cap
\mbox{\boldmath
$Q$}$. Let $\Delta$ be a {\bf Q}-boundary such that $K_X+\Delta$ is {\bf Q}-Cartier. 
Take a projective resolution $\mu:Y\rightarrow X$ 
such that each components of the support of 
$\mu^{-1}_{\ast}\Delta+\sum_{i\in I}E_i$ are smooth and cross normally, where $\{E_i\}_{i\in I}$ is a set of all the exceptional divisors of $\mu$
and put  
$d_i:=\mbox{mult}_{E_i}(K_Y+\mu^{-1}_{\ast}\Delta+\sum_{i\in I}E_i-\mu^{\ast}(K_X+\Delta))\in \mbox{\boldmath $Q$}$ for $i\in I$. 
The pair $(X,\Delta)$ is said to be {\it divisorially log terminal}, 
if all $d_i$ are positive and the exceptional loci of $\mu$ is purely one codimensional for some
$\mu$.  $(X,\Delta)$ is said to be {\it purely log terminal}, if all $d_i$ are positive for any
$\mu$. Let $\Delta:=\sum_{i}m_i\Gamma_i$ be the irreducible decomposition of $\Delta$.
$\lfloor\Delta\rfloor:=\sum_{m_i=1}\Gamma_i$ is called the {\it reduced part of $\Delta$} and 
$\{\Delta\}:=\sum_{m_i<1}d_i\Gamma_i$ is called the {\it fractional part of $\Delta$}. Let
$\Gamma$ be a normal prime divisor which does not contained in $\mbox{Supp }\Delta$.
$\mbox{Diff}_{\Gamma}(\Delta):=(K_X+\Gamma+\Delta)|_{\Gamma}-K_{\Gamma}$ is a {\bf Q}-divisor
called a {\it different} (see \cite{utah} and \cite{shokurov}).  
\\[2ex]
Let us fix our terminology from the category theory. By a {\it projective system}, we mean a category 
$\mathcal{I}$ such that 
$\mbox{Hom }_{\mathcal{I}}(\lambda,\mu)$ is empty or consists of exactly one element
$f_{\lambda,\mu}$ satisfying $f_{\lambda,\mu}\circ f_{\mu,\nu}=f_{\lambda,\nu}$ for any
$\lambda$, $\mu$,  $\nu\in \mbox{Ob }\mathcal{I}$. An object $\alpha\in \mbox{Ob }\mathcal{I}$ (
resp. $\omega\in \mbox{Ob }\mathcal{I}$ ) is called an {\it initial object} ( resp. a {\it final
object} ) if $\mbox{Card Hom }_\mathcal{I}(\alpha, \lambda)=1$ ( resp. $\mbox{Card Hom }_
\mathcal{I}(\lambda,\omega)=1$) for any $\lambda\in \mbox{Ob }\mathcal{I}$. A projective system $\mathcal{I}$
is said to be {\it cofilterd} if, for any given two objects $\lambda$, $\mu\in \mbox{Ob }\mathcal{I}$,
there exists $\nu\in \mbox{Ob }\mathcal{I}$ with $\mbox{Card Hom }_{\mathcal{I}}(\nu,\lambda)=\mbox{Card Hom }_{\mathcal{I}}(\nu,\mu)=1$.  A covariant functor $\Phi:\mathcal{I}^{\circ}\rightarrow
\mathcal{I}^{\prime \circ}$ between injective systems $\mathcal{I}^{\circ}$ and $\mathcal{I}^{\prime \circ}$ is said to be {\it cofinal}, if, for any given
$\lambda^{\prime}\in \mbox{Ob }\mathcal{I}^{\prime \circ}$, there exists $\lambda\in \mbox{Ob }\mathcal{I}^{\circ}$
such that $\mbox{Card Hom }_{\mathcal{I}^{\prime \circ}}(\lambda^{\prime},\Phi(\lambda))=1$. 
We shall also say that a projective subsystem $\mathcal{I}^{\prime}$ in a projective system $\mathcal{I}$ is cofinal in $\mathcal{I}$ 
if the dual embedding functor from $\mathcal{I}^{\prime \circ}$
to $\mathcal{I}^{\circ}$ is cofinal.  (see \cite{artin-mazur}, Appendix (1.5), \cite{sga4}, Expos\'e I, Definition 2.7 and Definition
8.1.1). 

\vskip 5mm

\noindent{\bf Acknowledgment.} The author would like to express his deep gratitude to to Prof.
Akira Fujiki for helping him to work on the analytic spaces and for giving him valuable
suggestions, to Prof. Makoto Namba for kindly showing him the beautiful book \cite{namba}, which enabled him to improve the first draft of this paper.

\section{Definition of $D$-local algebraic
fundamental groups} 
Let
$(X,p)$ be a germ of irreducible normal complex analytic spaces and let $\mathcal{M}_X$ denote the field of
germs of meromorphic functions on $X$. In what follows, we fix an algebraic
closure $\overline{\mathcal{M}_X}$ of $\mathcal{M}_X$ and the inclusion $i_X:\mathcal{M}_X\rightarrow \overline{\mathcal{M}_X}$. 
Take any $D\in \mbox{Div}_{\mbox{\boldmath $Q$}}X$ and fix it. Recall that a holomorphic map between complex analytic spaces
 is said to be {\it finite}, if it is proper
with discrete finite fibres.
\begin{df}{\em A finite
surjective morphism $f:Y\rightarrow X$, where $Y$ is a germ of irreducible normal complex analytic
spaces such that 
$f^{\ast}D$ is integral and Cartier, is called an {\it connected index one cover
with respect to }$D$. A connected index one cover $f:Y\rightarrow X$ with respect to $D$ is called {\it connected index one Galois cover
with respect to }$D$ if $f$ is Galois.
}
\end{df}
\begin{rk}{\em The above terminology is different from the one in \cite{shokurov}, in which \lq\lq an
index one cover'' means a canonical or log canonical cover in our terminology.
}\end{rk}
\begin{df}{\em Connected index one covers (resp. connected index one Galois covers) of $X$ with respect to $D$ form a full subcategory of complex analytic germs 
denoted by
$\mathcal{I}_{1}^m(X;D)$ (resp. $\mathcal{I}_{1}^mG(X;D)$). For $(Y,f)\in \mbox{Ob }\mathcal{I}_{1}^m(X;D)$, an injective homomorphism $i_Y:\mathcal{M}_Y\rightarrow
\overline{\mathcal{M}_X}$, where $\mathcal{M}_Y$ is the meromorphic
function field of $Y$, such that $i_Y\circ f^{\ast}=i_X$ is called a {\it pointing}. Triplet
$(Y, f, i_Y)$ composed of  $(Y,f)\in \mbox{Ob }\mathcal{I}_{1}^m(X;D)$ and a pointing $i_Y$ are called {\it pointed connected index one
covers with respect to }$D$. Pointed connected index one covers (resp. pointed connected index one Galois covers) with respect to $D$ and morphisms $f_{\lambda,\mu}\in \mbox{Hom
}_{\mathcal{I}_{1}^m(X;D)}((Y_{\mu},f_{\mu}),(Y_{\lambda},f_{\lambda}))$ satisfying
$i_{Y_{\mu}}\circ f_{\lambda,\mu}^{\ast}=i_{Y_{\lambda}}$ form a projective system denoted by $\mathcal{I}_{1}^m(X;D)^p$ (resp. $\mathcal{I}_{1}^mG(X;D)^p$). 
}\end{df}
Non-zero $\mbox{\boldmath $C$}$-algebra $\mathcal{A}$ is called
{\it a complex analytic ring} if there exists a surjective
$\mbox{\boldmath $C$}$-algebra homomorphism $\mathcal{O}^{\rm an}_{E,0}
\rightarrow \mathcal{A}$, where $E\simeq\mbox{\boldmath $C$}^n$ for some $n$. Recall that the category of complex analytic rings $\mathcal{A}$ which are finite $\mathcal{O}_X$-modules and that
category of germs of complex analytic spaces which are finite over $X$ are dual to each other via the contravariant functor $\mbox{Specan}_X$ (see
\cite{fischer} and \cite{grothendieck}, VI). The
structure morphism $f^{\ast}:\mathcal{O}_X\rightarrow \mathcal{A}$ is injective if and only if $f:\mbox{Specan}_X\mathcal{A}\rightarrow X$
is surjective by the Remmert's proper mapping theorem (see, for example, \cite{fischer}, 1.18). Let $\varphi$ be a meromorphic function on $X$
such that $\mathcal{O}_X(-rD)=\varphi\mathcal{O}_X$ and let 
$\pi:\tilde X\rightarrow X$ be the canonical cover with respect to $D$ obtained by
taking a
$r$-th root of $\varphi$, where $r:=\mbox{ind}_pD$ and fix a pointing $i_{\tilde X}$. Take any
$(Y,f,i_Y)\in
\mbox{Ob }\mathcal{I}_{1}^m(X;D)^p$. For simplicity, assume that 
$\mathcal{M}_X\subset\mathcal{M}_{\tilde X}\subset\overline{\mathcal{M}_X}$ and $\mathcal{M}_X\subset\mathcal{M}_Y\subset\overline{\mathcal{M}_X}$. The assumption on
$Y$ implies that there exists a meromorphic function $\psi$ on $Y$ such that
$\varphi\mathcal{O}_Y=\psi^{r}\mathcal{O}_Y$, which implies that there
exists a unit
$u\in \mathcal{O}_Y^{\times}$ such that $\psi^r=u\varphi$. Since 
$\mathcal{O}_Y$ is a henselian local ring whose residue field is the complex number field (see, for example, \cite{abhyankar}, Ch. III, \S 20,
Proposition 20.6), we see that $\root{r}\of{u}\in \mathcal{O}_Y^{\times}$, hence
$\mathcal{M}_X(\root{r}\of{\varphi})\subset \mathcal{M}_Y$. Consequently, there exists a $\mathcal{O}_X$-homomorphism 
$\pi_{\ast}\mathcal{O}_{\tilde
 X}\rightarrow f_{\ast}\mathcal{O}_{Y}$ which induces a finite surjective morphism
$\varpi_Y:Y=\mbox{Specan}_Xf_{\ast}\mathcal{O}_Y\rightarrow \tilde X=\mbox{Specan}_{X}\pi_{\ast}
\mathcal{O}_{\tilde X}$ satisfying $f=\pi\circ\varpi_Y$. The above argument implies that 
$\mbox{Card Hom }_{\mathcal{I}_{1}^m(X;D)^p}((Y,f,i_Y),(\tilde X,\pi,i_{\tilde X}))=1$ for any $(Y,f,i_Y)\in
\mbox{Ob }\mathcal{I}_{1}^m(X;D)^p$, that is, $(\tilde X,\pi,i_{\tilde X})$ is a final object, or equivalently, a colimit of 
$\mathcal{I}_{1}^m(X;D)^p$. Let $\varpi_Y(i_Y,i_{\tilde X})$ denote the element of 
$\mbox{Hom }_{\mathcal{I}_{1}^m(X;D)^p}((Y,f,i_Y),(\tilde X,\pi,i_{\tilde X}))$.

\begin{df}{\em A connected index one cover $f:Y\rightarrow X$ with respect to $D$ is called a {\it strict 
index one cover with respect to} $D$, if $\varpi_Y(i_Y,i_{\tilde X})$ is \'etale in codimension one for any pointings $i_Y$, $i_{\tilde X}$.
}\end{df}

\begin{rk}\label{rk:et1} {\em $\varpi_Y(i_Y,i_{\tilde X})$ is \'etale in codimension one if and only if $\varpi_Y(i_Y,i_{\tilde X})$ is \'etale over $\mbox{Reg }X$ by the
purity of branch loci (see 
\cite{abhyankar}, V, \S 39, (39.8) or \cite{fischer}, 4.2).
}\end{rk}

\begin{df}{\em A strict 
index one cover $f:Y\rightarrow X$ with respect to $D$ is called a {\it strict 
index one Galois cover} with respect to $D$, if $f$ is Galois.
}\end{df}

\begin{rk}{\em Take another pointings $i_Y^{\prime}$ and $i_{\tilde X}^{\prime}$ of
$(Y,f)$,
$(\tilde X,\pi)\in \mbox{Ob }\mathcal{I}_{1}^m(X;D)^p$ respectively and assume that $f:Y\rightarrow X$ is Galois. Then, by the Galois theory, 
there exist two 
isomorphisms $\alpha(i_Y,i_Y^{\prime})\in \mbox{Hom }_{\mathcal{I}_{1}^m(X;D)^p}((Y,f,i_Y), (Y,f,i_Y^{\prime}))$ and
$\beta(i_{\tilde X},i_{\tilde X}^{\prime})\in \mbox{Hom }_{\mathcal{I}_{1}^m(X;D)^p}((\tilde X,\pi,i_{\tilde X}),
(\tilde X,\pi,i_{\tilde X}^{\prime}))$ such that the following diagram in $\mathcal{I}_{1}^m(X;D)^p$ commutes.

\renewcommand{\normalbaselines}{\baselineskip20pt
\lieskip3pt \lineskiplimit3pt}
\newcommand{\mapright}[1]{\smash{\mathop{
\hbox to 1cm{\rightarrowfill}}\limits^{#1}}}

\newcommand{\mapleft}[1]{\smash{\mathop{
\hbox to 1cm{\leftarrowfill}}\limits_{#1}}}

\newcommand{\mapdown}[1]{\Big\downarrow
\llap{$\vcenter{\hbox{$\scriptstyle#1\,$}}$ }}

\newcommand{\mapdownr}[1]{\Big\downarrow
\rlap{$\vcenter{\hbox{$\scriptstyle#1\,$}}$ }}

\newcommand{\mapup}[1]{\Big\uparrow\rlap{$\vcenter{\hbox{$\scriptstyle#1$}}$ }}
\[\begin{array}{ccc}
(Y, f,i_Y)&\mapright{\alpha(i_Y,i_Y^{\prime})}&(Y, f,i_Y^{\prime})\\
\mapdown{\varpi_Y(i_Y,i_{\tilde X})}&&\mapdownr{\varpi_Y(i_Y^{\prime},i_{\tilde X}^{\prime})} \\
(\tilde X,\pi,i_{\tilde X})&\mapright{\beta(i_{\tilde X},i_{\tilde X}^{\prime})} & (\tilde X,\pi,i_{\tilde X}^{\prime})
\end{array}\]

Therefore, $(Y,f)\in \mbox{Ob }\mathcal{I}_{1}^m(X;D)^p$ is a strict 
index one Galois cover if one of $\varpi_Y(i_Y,i_{\tilde X})\in \mbox{Hom }_{\mathcal{I}_{1}^m(X;D)}((Y,f),(\tilde X,\pi))$ is \'etale in codimension one. One
can also check easily that the same holds even if $f$ is not Galois.
}
\end{rk}

Let $\mathcal{I}_1^{\dag}(X;D)$ (resp, $\mathcal{I}_1^{\dag}(X;D)^p$) denote the full subcategory of $\mathcal{I}_{1}^m(X;D)$ (resp.
projective subsystem $\mathcal{I}_{1}^m(X;D)^p$) whose objects are strict index one covers with respect to
$D$ (resp. pointed strict index one covers with respect to
$D$) and let $\mathcal{I}_1^{\dag}G(X;D)$ (resp, $\mathcal{I}_1^{\dag}G(X;D)^p$) denote the full subcategory of $\mathcal{I}_1^{\dag}(X;D)$ (resp.
projective subsystem $\mathcal{I}_1^{\dag}(X;D)^p$) whose objects are strict index one Galois covers with respect to
$D$ (resp. pointed strict index one Galois covers with respect to
$D$).  Let $f_{\lambda,\mu}:(Y_{\mu},f_{\mu}, i_{Y_{\mu}})\rightarrow (Y_{\lambda},f_{\lambda},i_{Y_{\lambda}})$ be a morphism
in
$\mathcal{I}_1^{\dag}G(X;D)^p$ and assume that 
$\mathcal{M}_X\subset\mathcal{M}_{Y_{\lambda}}\subset \mathcal{M}_{Y_{\mu}}\subset 
\overline{\mathcal{M}_X}$ for simplicity. Then by the Galois theory, there
exists a canonical surjective homomorphism
$g_{\lambda,\mu}:\mbox{Gal }(\mathcal{M}_{Y_{\mu}}/\mathcal{M}_X)\rightarrow \mbox{Gal }(\mathcal{M}_{Y_{\lambda}}/\mathcal{M}_X)$ which is
nothing but the restriction map. Thus Galois groups $\mbox{Gal }(Y/X):=\mbox{Gal }(\mathcal{M}_Y/\mathcal{M}_X)$ for $(Y,f,i_Y)\in \mathcal{I}_1^{\dag}G(X;D)^p$ 
form a projective system with the induced morphisms from $\mathcal{I}_1^{\dag}G(X;D)^p$.  

\begin{df}\label{df:log pi1}{\em  We define a profinite group $\hat\pi_{1,X,p}^{{\rm loc}}[D]$
by 
$$
\hat\pi_{1,X,p}^{{\rm loc}}[D]:=\mbox{projlim}_{(Y,f,i_Y)\in {\rm Ob\ }\mathcal{I}_1^{\dag}G(X;D)^p}\mbox{Gal }(Y/X) 
$$
which is called {\it a divisorial local algebraic fundamental group} with respect to $D$, or {\it a $D$-local algebraic fundamental group}.
}\end{df}
\begin{rk} Obviously, $D$-local algebraic fundamental groups do not depend on the choice of inclusions $i_X:\mathcal{M}_X\rightarrow \overline{\mathcal{M}_X}$.

\end{rk}
\begin{rk}{\em Our profinte groups can not be defined directly using Grothendieck's theory \lq\lq
cat\'egories galoisiennes'' because we can not find a suitable category for our theory. The
problem is the existence of a final and an initial object as in the axioms (G 1) and (G 2)  (see
\cite{sga1}, Expos\'e V, \S 4).
}\end{rk}

\begin{rk}{\em Take two {\bf Q}-Cartier, {\bf Q}-divisors $D_1$ and $D_2$ such that $D_1\sim D_2$,
that is, $D_1-D_2$ is integral and Cartier. Then it is obvious by the definition that $\mathcal{I}_1^{m (\dag)}(X;D_1)=
\mathcal{I}_1^{m (\dag)}(X;D_2)$ and  
$\hat\pi_{1,X,p}^{{\rm loc}}[D_1]=\hat\pi_{1,X,p}^{{\rm loc}}[D_2]$. In other words, $\hat\pi_{1,X,p}^{{\rm loc}}[D]$ depends only on the
class $[D]\in \mbox{Div}_{\mbox{\boldmath $Q$}}X/\mbox{Div }X$. }
\end{rk}

\section{Comparison theorem}
In this section, 
we shall compare $D$-local algebraic fundamental groups with certain topological invariants. For a germ of normal complex analytic spaces $(X,p)$, put 
$\mbox{Reg }X:=\mbox{projlim}_{p\in \mathcal{U};{\rm open }}\mbox{Reg }\mathcal{U}$, where $\mbox{Reg }\mathcal{U}$ is the smooth loci of $\mathcal{U}$ and  
let $\pi_1^{{\rm loc}}(\mbox{Reg }X)$ denote the local fundamental group defined as $\mbox{projlim}_{p\in \mathcal{U};{\rm open }}\pi_1(\mbox{Reg }\mathcal{U})$ 
and $\hat\pi_1^{{\rm loc}}(\mbox{Reg }X)$ its profinite completion which is called the {\it local algebraic fundamental group} of $(X,p)$. Let $\Sigma$ be an analytically closed
proper subset of $X$.  According to Prill (\cite{prill}, \S IIB), there exists a contractible open neighbourhood $U$ of $p$ such that there exists a neighbourhood basis
$\{U_{\lambda}\}_{\lambda\in
\Lambda}$ of $p$ satisfying the condition that
$U_{\lambda}\setminus\Sigma$ is a deformation retract of $U\setminus \Sigma$ for any $\lambda\in
\Lambda$. By the definition, we have $\pi_1(U\setminus\Sigma)=\mbox{projlim}_{p\in \mathcal{U};{\rm open }}\pi_1(\mathcal{U}\setminus\Sigma)$. 
We call such $U$ as above a {\it Prill's good neighbourhood with regard to
$\Sigma$} and we say that 
$\{U_{\lambda}\}_{\lambda\in \Lambda}$ is a {\it neighbourhood basis associated with $U$}. Recall that $U_{\lambda}$ is also a Prill's good neighbourhood with
regard to $\Sigma$ and for any two Prill's good neighbourhood $U$ and $U^{\prime}$, $U\setminus\Sigma$ and $U^{\prime}\setminus\Sigma$ have the same homotopy type.
In particular, we have $\pi_1^{{\rm loc}}(\mbox{Reg }X)\simeq\pi_1(\mbox{Reg }U)$ for a Prill's good neighbourhood $U$ with
regard to $\mbox{Sing }X$. To introduce the generalized notion of local fundamental groups, let us briefly review here the theory of universal ramified coverings due to M. Kato
(\cite{kato}), M. Namba (\cite{namba}) and J.P. Serre (\cite{serre}, Appendix 6.4) according to M. Namba.  Let $B$ be an integral effective divisor on a connected complex manifold $M$
and let
$B:=\sum_{i\in I}b_iB_i$ be the irreducible decomposition of
$B$.  Fix a base point $x\in M\setminus\mbox{Supp }B$ and let $\gamma_i$ be a loop which starts from $x$ and goes around $B_i$ once in a counterclockwise direction with the center
being a smooth point of
$\mbox{Supp }B$ on $B_i$. Let $\mathcal{N}(M,B,x)\subset \pi_1(M\setminus \mbox{Supp }B,x)$ denote the normal subgroup generated by all the conjugates of the loops
$\{\gamma_i^{b_i}\}_{i\in I}$. Recall that $\mathcal{N}(M,B,x)$ is known to be independent from the choice of such loops. We define a {\it $B$-fundamental group} of $M$ by putting  
$$
\overline\pi_1^B(M,x):=\pi_1(M\setminus \mbox{Supp }B,x)/\mathcal{N}(M,B,x).
$$
A finite covering $f:N\rightarrow M$ from a connected normal complex analytic space $N$ which is \'etale over $M\setminus \mbox{Supp }B$ is said to be branching at most (resp.
branching ) at
$B$, if the ramification index $e_{\tilde B_{i,j}}(f)$ of $f$ at any prime divisor $\tilde B_{i,j}$ such that $f(\tilde B_{i,j})=B_i$ divides (resp. is equals to) $b_i$ for any $i\in
I$. Let
$FC^{\leq B}(M)$ (resp. $FC^B(M)$ ) denote the category of finite coverings over $M$  branching at most (resp. branching) at $B$.  Let $FGC^{\leq B}(M)$ (resp. $FGC^B(M))$ denote the
full subcategory of $FC^{\leq B}(M)$ whose objects consists of Galois covers over $M$. Triplet $(N, f, y)$, where $(N,f)\in \mbox{Ob }FC^{\leq B}(M)$ and $y\in f^{-1}(x)$ are called
{\it pointed finite coverings branching at most at $B$}. Pointed finite coverings branching at most at $B$ and morphisms $f_{\lambda,\mu}\in
\mbox{Hom}_{FC^{\leq B}(M)}((N_{\mu},f_{\mu}),(N_{\lambda},f_{\lambda}))$ such that
$f_{\lambda,\mu}(y_{\mu})=y_{\lambda}$, where $(N_{\mu},f_{\mu}, y_{\mu})$ and $(N_{\lambda},f_{\lambda}, y_{\lambda})$ are two pointed finite coverings branching at most at $B$ 
form a
projective system denoted by $FC^{\leq B}(M)^p$. We also define the projective subsystems $F(G)C^{(\leq) B}(M)^p$ in the same way. 
From \cite{namba}, Lemma 1.3.1, Theorem 1.3.8 and 
Theorem 1.3.9, we see that there exists a canonical functor $\Psi$ from $FC^{\leq B}(M)^p$ to the projective system of subgroups of finite indices in $\overline\pi_1^B(M,x)$ such that 
$$
\Psi((N,f,y))=f_{\ast}\pi_1(N\setminus \mbox{Supp }f^{-1}B,y)/\mathcal{N}(M,B,x)\subset \overline\pi_1^B(M,x)
$$
for 
$(N,f,y)\in \mbox{Ob}FC^{\leq B}(M)^p$ and that the functor $\Psi$ defines an equivalence between the above two projective systems. Thus by using the basic group theory, we obtain the
following lemma. 

\begin{lm}\label{lm:fgc}
$FGC^{\leq B}(M)^p$ $($ resp. $FGC^B(M)^p$ $)$ is cofilterd and cofinal in $FC^{\leq B}(M)^p$ $($ resp. $FC^B(M)^p$ $)$ and hence in particular, $FGC^B(M)^p$ is cofinal in
$FGC^{\leq B}(M)^p$ if
$FGC^B(M)^p$ is not empty. 
\end{lm}

\begin{rk}{\em Let $\overline\pi_1^B(M,x)^{\wedge}$ denote the profinite completion of $\overline\pi_1^B(M,x)$ called the {\it $B$-algebraic fundamental group} of $M$.
Assume that $FGC^B(M)^p$ is not empty. Then by Lemma~\ref{lm:fgc}, we have
$$
\overline\pi_1^B(M,x)^{\wedge}\simeq \mbox{projlim}_{(N,f,y)\in {\rm Ob}FGC^B(M)^p}\mbox{Gal }(N/M), 
$$
where $\mbox{Gal }(N/M):=\overline\pi_1^B(M,x)/f_{\ast}\overline\pi_1^B(N,y)$.
}\end{rk}

\begin{df}{\em 
For a germ of normal complex analytic spaces $(X,p)$ and $B\in \mbox{Weil }X$, we 
define a $B$-local fundamental group of $X$ with respect to $B$ as follows:
$$
\overline\pi_{1,{\rm loc}}^B(\mbox{Reg }X):=\mbox{projlim}_{p\in \mathcal{U};{\rm open}}\overline\pi_1^B(\mbox{Reg }\mathcal{U}). 
$$
Moreover, by $\overline\pi_{1,{\rm loc}}^B(\mbox{Reg }X)^{\wedge}$, we
mean the profinite completion of $\overline\pi_{1,{\rm loc}}^B(\mbox{Reg }X)$.
}\end{df}

\begin{rk}{\em We note that apparently, we have $\overline\pi_{1,{\rm loc}}^B(\mbox{Reg }X)=\pi_1^{{\rm loc}}(\mbox{Reg }X)$ and $\overline\pi_{1,{\rm loc}}^B(\mbox{Reg
}X)^{\wedge}=\hat\pi_1^{{\rm loc}}(\mbox{Reg }X)$ if $B=0$.
}\end{rk}

Before arguing about the comparison with $D$-local algebraic fundamental groups and $B$-local algebraic fundamental groups, we introduce certain categories containing the categories of
connected or strict index one covers as a full subcategories as an intermediate step.

\begin{df}{\em  A finite
surjective morphism $f:Y\rightarrow X$, where $Y$ is a germ of irreducible normal complex analytic
spaces such that 
$f^{\ast}D\in \mbox{Weil }Y$ is called a {\it connected integral cover
with respect to }$D$. A connected integral cover $f:Y\rightarrow X$ with respect to $D$ is called a {\it strict integral
cover},  if $e_{\tilde \Gamma}(f)=e_{\Gamma}(D)$ for any prime divisors $\tilde\Gamma$ on $Y$ and $\Gamma$ on $X$ such that $f(\tilde\Gamma)=\Gamma$, where $e_{\tilde \Gamma}(f)$
denotes the ramification index of $f$ at $\tilde\Gamma$ and $e_{\Gamma}(D):=[\mbox{\boldmath $Z$}(\mbox{mult}_{\Gamma}D):\mbox{\boldmath
$Z$}(\mbox{mult}_{\Gamma}D)\cap \mbox{\boldmath $Z$}]\in \mbox{\boldmath $N$}$. By $\mbox{Int}^m(X;D)$ (
resp.
$\mbox{Int}^{\dag}(X;D)$ ), we mean a category of connected integral covers ( resp. a category of strict integral covers ) with respect to $D$. We shall also define categories
$\mbox{Int}^{m(\dag)}(G)(X;D)^{(p)}$ similarly as before. }\end{df}

\begin{rk}{\em By the definition, one can see that the category of strict integral covers 
$\mbox{Int}^{\dag}(X;D)$ contains the category of strict index one covers $\mathcal{I}_1^{\dag}(X;D)$ as a full subcategory.  Note that in particular, we have $\mbox{Int}^{\dag}(X;0)=\mathcal{I}_1^{\dag}(X;0)$.
}\end{rk}

Let $\mathcal{X}$ be an arcwise connected, locally arcwise
connected, Hausdorff topological space, A continuous map $f:\mathcal{Y}\rightarrow \mathcal{X}$ from a Hausdorff topological space $\mathcal{Y}$ with discrete finite fibres
is called a {\it finite topological covering} if for any $x\in\mathcal{X}$, there exists an arcwise connected open neighbourhood $\mathcal{U}$ of
$x\in\mathcal{X}$ such that the restriction of $\pi$ to each arcwise connected component of of $\pi^{-1}(\mathcal{U})$ gives a homeomorphism
onto $\mathcal{U}$. The following lemma is nothing but a consequence from the first covering homotopy theorem (see \cite{steenrod}, 11.3).

\begin{lm}\label{lm:lifting of homotopy} Let $f:\mathcal{Y}\rightarrow \mathcal{X}$ be a connected finite topological covering. Assume that
$\mathcal{X}$ is paracompact and let
$\mathcal{Z}\subset \mathcal{X}$ be a topological subspace which is a deformation retract of $\mathcal{X}$. Then 
$\tilde{\mathcal{Z}}:=\pi^{-1}(\mathcal{Z})\subset \mathcal{Y}$ is a deformation retract of $\mathcal{Y}$. In particular, 
$\tilde{\mathcal{Z}}$ is also 
arcwise connected and
$\pi_1(\tilde{\mathcal{Z}})=\pi_1(\mathcal{Y})$.
\end{lm}

For a germ of normal complex analytic space $(X,p)$, let $U$ be a Prill's good
neighbourhood with regard to a proper analytically closed subset $\Sigma\subset X$ and let $\{U_{\lambda}\}_{\lambda\in\Lambda}$ be a neighbourhood basis
associated with $U$. We put $U^{-}:=U\setminus \Sigma$ and $U^{-}_{\lambda}:=U_{\lambda}\setminus \Sigma$. Take any connected finite topological covering $f^{-}:V^{-}\rightarrow
U^{-}$.  Then $V^-$ has the unique analytic structure such
that $f^-:V^-\rightarrow U^{-}$ is \'etale. By the Grauert-Remmert's theorem, $f^-$
extends uniquely to a finite cover
$f:V\rightarrow U$, where $V$ is a normal complex analytic space such that $f^{-1}(U^{-})=V^-$ 
(see \cite{grauert-remmert}, \S 2, Satz 8 and \cite{sga1}, XII,
Theorem 5.4). Recall here that $f^{-1}(p)$ consists of exactly one point, for if $f^{-1}(p)=\{q_1,\dots ,q_n\}$ and $n\geq
2$, where
$q_i$ are distinct from each other, then by \cite{fischer}, 1.10, Lemma 2, there exists $\lambda\in \Lambda$ such that
$f^{-1}(U_{\lambda})=\coprod_{i=1}^n W_{\lambda,i}$, where $W_{\lambda,i}$ is an open neighbourhood of $q_i$ for $i=1,\dots, n$.
Since $f^{-1}(U_{\lambda}^{-})\subset f^{-1}(U_{\lambda})$ is connected by
Lemma~\ref{lm:lifting of homotopy}, there exists $i$, say $i_0$, such that 
$f^{-1}(U_{\lambda}^{-})\subset W_{\lambda,i_0}$ and $f^{-1}(U_{\lambda}^{-})\cap W_{\lambda, i}$ is
empty if $i\neq i_0$, but which is absurd for  
$f^{-1}(U_{\lambda}^{-})\cap W_{\lambda, i}=W_{i,\lambda}\setminus f^{-1}(\Sigma)$ is non-empty for any $i$. Thus we see that any connected finite topological covering
$f^{-}:V^{-}\rightarrow U^{-}$ determines a finite surjective morphism $f:Y:=(V,q)\rightarrow X$ from a germ of normal complex analytic spaces $Y$ uniquely up to isomorphisms, where
$f^{-1}(p)=\{q\}$. For two connected finite topological coverings $f^{-}_1:V^{-}_1\rightarrow U^{-}$ and $f^{-}_2:V^{-}_2\rightarrow U^{-}$, 
let $f_i:V_i\rightarrow U^{-}$ be the
extended finite covers and $f:Y_i\rightarrow X$ be the corresponding finite surjective morphisms as above for
$i=1$, $2$. By
\cite{sga1}, Expos\'e XII, Proposition 5.3, the restriction map 
$r:\mbox{Hom}_U(V_1,V_2)\rightarrow \mbox{Hom}_{U^{-}}( V^-_1, V^-_2 )$ is
bijective and composed with the canonical injection $\mbox{Hom}_U(V_1,V_2)\rightarrow \mbox{Hom}_X(Y_1,Y_2)$, $r^{-1}$ gives a canonical injection $\mbox{Hom}_{U^{-}}( V^-_1, V^-_2)\rightarrow  \mbox{Hom}_{X}(Y_1,Y_2)$. 
From the above argument, we see that we have a canonical faithful functor $\mathcal{P}$ called a {\it Prill functor} from the category of connected topological finite coverings over $U^{-}$ denoted by $FT(U^{-})$ to the category of germs of normal complex analytic spaces which is finite over $X$ and \'etale outside over $\Sigma$ denoted by $FC(X,\Sigma)$.
\begin{lm}\label{lm:prill-fun} A Prill functor defines an equivalence between the categories $FT(U^{-})$ and $FC(X,\Sigma)$.
\end{lm}
{\it Proof. } We note that the restriction functor $\mathcal{R}_{\lambda}:FT(U^{-})\rightarrow FT(U_{\lambda}^{-})$ defines an equivalence of categories between
$FT(U^{-})$ and
$FT(U_{\lambda}^{-})$ since these are known to be determined up to equivalences by the corresponding fundamental
groups. Put 
$(V_{i,\lambda}^-, f_{i,\lambda}^-):=\mathcal{R}_{\lambda}((V_i^-,f_i^-))\in \mbox{Ob }FT(U_{\lambda})$ and
$V_{i,\lambda}:=f^{-1}_i(U_{\lambda})$ for $i=1,2$. Note also that the canonical map 
$$
\mbox{injlim}_{\lambda\in \Lambda} \mbox{Hom
}_{U_{\lambda}}(V_{1,\lambda},V_{2,\lambda})\rightarrow \mbox{Hom }_{FC(X,\Sigma)}((Y_1,f_1),(Y_2,f_2))
$$ is bijective and 
 that we have 
$$
\mbox{injlim}_{\lambda\in \Lambda} \mbox{Hom }_{FT(U^{-}_{\lambda})}((V_{1,\lambda}^-,f_{1,\lambda}^-),(V_{2,\lambda}^-,f_{2,\lambda}^-))=\mbox{Hom
}_{FT(U^{-})}((V_1^-,f_1^-),(V_2^-,f_2^-)).
$$ Therefore, we conclude that the canonical map 
$$\mbox{Hom }_{FT(U^{-})}( (V^-_1,f^-_1), (V^-_2,f^-_2) )\rightarrow  \mbox{Hom }_{FC(X,\Sigma)}((Y_1,f_1),(Y_2,f_2))
$$ is
bijective, which implies that the functor $\mathcal{P}$ is faithfully full. Take any $(Y,f)\in \mbox{Ob }FC(X,\Sigma)$. Then $f$ is
represented by a finite cover $f:V_{\lambda}\rightarrow U_{\lambda}$ for some $\lambda \in \Lambda$, where $V_{\lambda}$ is
connected. Since 
$f$ is
\'etale over
$U_{\lambda}^-$ and $V_{\lambda}^-:=f^{-1}(U_{\lambda}^-)$ 
is also connected, we obtain an object 
$(V_{\lambda}^-, f|_{V_{\lambda}^-})\in \mbox{Ob }FT(U^{-}_{\lambda})$ which goes to $(Y,f)\in \mbox{Ob }FC(X,\Sigma)$ via 
$\mathcal{P}\circ \mathcal{R}_{\lambda}^{-1}$.  Thus we conclude that $\mathcal{P}$ is essentially surjective\footnote{This is the English translation of \lq\lq essentiellement surjectif'' in French, and has the
same meaning as 
\lq\lq representative'' in \cite{mitchell}, II, 4.}, and hence
$\mathcal{P}$ defines an equivalence.
\hfill\bsquare
\begin{rk}\label{rk:prill-fun} {\em It is obvious that a Prill functor also defines an equivalence between the full
subcategory of Galois objects of $FT(U^{-})$ and $FC(X,\Sigma)$.  Note that giving a pointing to an object of $FT(U^{-})$ and $FC(X,\Sigma)$ has essentially the same meaning since the
number of pointings for $(V^-,f^-)\in \mbox{Ob }FT(U^{-})$  and the number of pointings for $(Y,f)\in \mbox{Ob }FC(X,\Sigma)$ are both $\deg f=\deg f^-$.  Thus we see
that $\mathcal{P}$ induces an equivalence between 
$FT(U^{-})^p$ and $FC(X,\Sigma)^p$. 
}\end{rk}
\begin{notation}\label{notation:Dvee}{\em
Let $\mathcal{B}_X(D)$ be the set of all the prime divisors on $X$ such that $e_{\Gamma}(D)>1$ and define the Weil divisor $D^{\vee}$ on $X$ by putting
$$D^{\vee}:=\sum_{\Gamma\in \mathcal{B}_X(D)}e_{\Gamma}(D)\Gamma.$$
}
\end{notation}
Combined with Lemma~\ref{lm:fgc}, Lemma~\ref{lm:prill-fun} yields the following proposition.
\begin{pr}\label{pr:top-int} There exists a canonical functor $\mathcal{P}:FC^{D^{\vee}}(\mbox{\em Reg }U)\rightarrow \mbox{\em Int}^{\dag}(X;D)$ which defines an equivalence between the
categories
$FC^{D^{\vee}}(\mbox{\em Reg }U)$ and $\mbox{\em Int}^{\dag}(X;D)$. In particular, $\mbox{\em Int}^{\dag}G(X;D)^p$
is cofilterd and cofinal in $\mbox{\em Int}^{\dag}(X;D)^p$. 
\end{pr}
\begin{rk}\label{rk:top-int}{\em From Proposition~\ref{pr:top-int}, we deduce that $\mbox{projlim}_{(Y,f,i_Y)\in {\rm Ob\ Int}^{\dag}G(X;D)^p}\mbox{Gal
}(Y/X)$ is isomorphic to $\overline\pi_{1,{\rm loc}}^{D^{\vee}}(\mbox{Reg }X)^{\wedge}$ in the category of profinite groups.
}\end{rk}
Let $\mathcal{C}$ be a category and let 
$X\in \mbox{Ob }\mathcal{C}$ be an object of $\mathcal{C}$ and $G\subset \mbox{Aut }X$ be
a subgroup of the automorphism group of $X$.
\begin{df}{\em An epimorphism $f:X\rightarrow Y$ in $\mathcal{C}$ is said to be {\it Galois with the Galois group} $G$,
if
$G=\mbox{Aut }_Y X:=\{\sigma\in \mbox{Aut }X|f\circ\sigma=f \}$ and for any morphism $f^{\prime}:X\rightarrow Y^{\prime}$ such
that
$G\subset \mbox{Aut }_{Y^{\prime}} X$, there exists a unique morphism $\varphi:Y\rightarrow Y^{\prime}$ satisfying
$f^{\prime}=\varphi\circ f$.
}
\end{df}
\begin{rk}{\em Assume that two Galois morphisms $f:X\rightarrow Y$ and $f^{\prime}:X\rightarrow Y^{\prime}$ with the Galois group
$G$ are given. Then by the universal mapping property, there exists an isomorphism $\varphi:Y\rightarrow Y^{\prime}$ such that
$f^{\prime}=\varphi\circ f$, that is, Galois morphisms with the Galois group
$G$ is unique up to this equivalence.
}\end{rk}
\begin{ex}{\em Let $\mathcal{F}:=(\mbox{Fields})$ be a category of fields such that $\mbox{Hom}_\mathcal{F}(K_1,K_2)$ is empty or
consists of inclusions for any $K_1$, $K_2\in \mbox{Ob }\mathcal{F}$. For any finite extension $i:K_1\rightarrow K_2$, 
$i$ is a Galois extension if and only if its dual $i^{\circ}:K_2^{\circ}\rightarrow K_1^{\circ}$ in the dual category 
$\mathcal{F}^{\circ}$ is Galois by the Galois theory.
}\end{ex}
\begin{df}{\em For any two morphisms $f\in \mbox{Hom}_{\mathcal{C}}(X,Y)$ and $g\in \mbox{Hom}_{\mathcal{C}}(Y,Z)$, we define a subgroup
$\mbox{Aut}^f_Z X\subset\mbox{Aut }X\times \mbox{Aut}_Z Y$ as $\mbox{Aut}^f_Z X:=\{(\tilde \sigma,\sigma)|f\circ\tilde
\sigma=\sigma\circ f\}$.
}\end{df}
The following lemma grew out of Prof. A. Fujiki's suggestion.
\begin{lm}\label{lm:galois} Let $f\in \mbox{\em Hom}_{\mathcal{C}}(X,Y)$ and $g\in \mbox{\em Hom}_{\mathcal{C}}(Y,Z)$ be two Galois
morphisms and assume that the second projection $p_2:\mbox{\em Aut}^f_Z X\rightarrow \mbox{\em Aut}_Z Y$ is surjective. Then
$h:=g\circ f$ is also Galois.
\end{lm}
{\it Proof.} Take any $h^{\prime}\in \mbox{Hom }_{\mathcal{C}}(X,Z^{\prime})$ with $\mbox{Aut}_Z X\subset \mbox{Aut}_{Z^{\prime}}X$.
Since $\mbox{Aut}_Y X\subset \mbox{Aut}_Z X$ and $f$ is Galois, there exists a morphism $\psi:Y\rightarrow Z^{\prime}$ such that
$\psi\circ f=h^{\prime}$. Take any $\sigma\in \mbox{Aut}_Z Y$. Then there exists $\tilde \sigma\in \mbox{Aut }X$ such that
$f\circ\tilde\sigma=\sigma\circ f$ by the assumption. Since $\tilde\sigma\in \mbox{Aut}_Z X\subset \mbox{Aut}_{Z^{\prime}}X$, we
have
$h^{\prime}=h^{\prime}\circ\tilde\sigma=\psi\circ f\circ \tilde\sigma=\psi\circ\sigma \circ f$, hence 
$\psi\circ f=\psi\circ\sigma
\circ f$. Since $f$ is an epimorphism, we deduce that $\psi=\psi\circ \sigma$, that is, $\mbox{Aut}_Z Y\subset
\mbox{Aut}_{Z^{\prime}}Y$. Thus we conclude that there exists a morphism $\varphi:Z\rightarrow Z^{\prime}$ such that $\varphi\circ
g=\psi$.  Obviously $\varphi$ satisfies $\varphi\circ h=h^{\prime}$. As for the uniqueness of $\varphi$, Let
$\varphi^{\prime}:Z\rightarrow Z^{\prime}$ be another morphism satisfying $\varphi^{\prime}\circ h=h^{\prime}$. Then $\psi\circ
f=h^{\prime}=\varphi^{\prime}\circ h=\varphi^{\prime}\circ g\circ f$, hence $\varphi^{\prime}\circ g=\psi=\varphi\circ g$. Since $g$
is also an epimorphism, we obtain $\varphi^{\prime}=\varphi$.
\hfill \bsquare
\begin{rk}\label{rk:lift}{\em The assumption in Lemma~\ref{lm:galois} is satisfied in the following two theoretically important cases.
\begin{enumerate}
\item Let $f$ and $g$ are finite Galois covers between normal algebraic varieties over an algebraically closed field or normal connected complex
analytic spaces. 
Assume that there exits a Zariski closed subset or an analytic subset
$\Sigma$ on
$Y$ with
$\mbox{codim}_Y \Sigma\leq 2$ such that the restriction $f^-$ of $f$ to $X^-:=X\setminus
f^{-1}(\Sigma)$ gives the algebraic universal cover of $Y^-:=Y\setminus \Sigma$, that is,
$\hat\pi_1(X^-)=\{1\}$. Moreover assume that $Y^-$ is invariant under the action of
$\mbox{Gal }(Y/Z)$. Take
any
$\sigma\in
\mbox{Gal }(Y/Z)$. Since $\sigma$ acts on
$Y^-$, there exists an automorphism
$\tilde\sigma^-$ on
$X^-$ such that $f^-\circ\tilde \sigma^-=\sigma\circ f^-$ by the property of algebraic universal
cover.
$\tilde\sigma^-$ extends uniquely to an automorphism $\tilde\sigma$ on
$X$ satisfying $f\circ\tilde \sigma=\sigma\circ f$ by the normality (see also \cite{catanese}, \S1 and \cite{gabi}, Lemma
2.1). 
\item Let $f$ and $g$ are finite Galois covers between germs of normal complex
analytic spaces. Assume that $X$ is obtained from
$Y$ by taking a $r$-th root of a primitive principal divisor $P=\mbox{div }\varphi$ on $Y$ such that 
$\mathcal{O}_Y(P)\subset\mathcal{M}_Y$ is invariant under the action of $\mbox{Gal }(Y/Z)$ (for the definition
of primitive principal divisors, see \cite{shokurov}, 2.3). Take any $\sigma\in \mbox{Gal }(Y/Z)$. Then by the assumption,
we have $\sigma ^{\ast}\varphi=u\varphi$ for some unit $u\in \mathcal{O}_Y^{\times}$. 
As in the previous argument, there exists a unit $v\in \mathcal{O}_Y^{\times}$ such that $v^r=u$. Since we can write $\mathcal{M}_X=\mathcal{M}_Y[T]/(T^r-\varphi)$, it is obvious that 
$\sigma
^{\ast}$ lifts to an automorphism $\tilde\sigma^{\ast}$ on $\mathcal{M}_X$ by putting $\tilde\sigma^{\ast}T=vT$.  Thus we see that any elements of $\mbox{Gal }(Y/Z)$ lift to 
elements of $\mbox{Gal }(X/Z)$. 
\end{enumerate} 
}\end{rk}
\begin{lm}[c.f., \cite{seidenberg2}]\label{lm:normalization} Let $\mathcal{A}$ be an integral complex analytic ring and $\mathcal{M}$ be its quotient field. Let
$\mathcal{A}_{\mathcal{L}}$ be the normalization of $\mathcal{A}$ in a finite extension field $\mathcal{L}$ of $\mathcal{M}$. Then 
$\mathcal{A}_{\mathcal{L}}$ is also an integral complex analytic ring which is a finite $\mathcal{A}$-module. 
\end{lm}
{\it Proof.} Recall that $\mathcal{A}$ is noetherian (\cite{grothendieck}, II, Proposition 2.3). \cite{seidenberg2}, Theorem 4
says that
$\mathcal{A}$ is N-1, hence N-2 by \cite{matsumura}, Ch. 12, Corollary 1, that is, 
$\mathcal{A}_{\mathcal{L}}$ is a finite $\mathcal{A}$-module. By \cite{seidenberg2}, Theorem 1, 
$\mathcal{A}$ is a finite $\mathcal{O}^{\rm an}_{C^n,0}$-module for some $\mbox{\boldmath $C$}^n$, hence so is 
$\mathcal{A}_{\mathcal{L}}$. Thus by \cite{seidenberg2}, Theorem 3, we conclude that $\mathcal{A}_{\mathcal{L}}$ is a complex analytic ring.
\hfill\bsquare
\begin{rk}\label{rk:normalization} {\em Lemma~\ref{lm:normalization} implies that if we are given a finite extension field $\mathcal{L}$ of the
meromorphic function field $\mathcal{M}_X$ of an irreducible germ of complex analytic spaces $X$, there exists a germ of normal
complex analytic spaces
$Y$ with a finite surjective morphism
$f:Y\rightarrow X$ such that $\mathcal{M}_Y$ is isomorphic to $\mathcal{L}$ over $\mathcal{M}_X$ and such $Y$ as above is uniquely determined up to isomorphisms over $X$.
}\end{rk}
The following proposition can be also derived from Proposition~\ref{pr:top-int}, but we shall give an algebraic proof for further research such as extending our theory to the positive
characteristic case.
\begin{pr}\label{pr:cofilterd} $\mathcal{I}_1^{\dag}G(X;D)^p$ is cofilterd and is cofinal in $\mathcal{I}_1^{\dag}(X;D)^p$ for any $D\in \mbox{\em Div}_{\mbox{\boldmath $Q$}}X$.
\end{pr}
{\it Proof. }Firstly, we prove the first statement. Take any two objects $(Y_{\lambda},f_{\lambda},i_{Y_{\lambda}})$, $(Y_{\mu},f_{\mu},i_{Y_{\mu}})\in \mathcal{I}_1^{\dag}G(X;D)^p$. 
Let
$\mathcal{L}:=i_{Y_{\lambda}}(\mathcal{M}_{Y_{\lambda}})
\lor i_{Y_{\mu}}(\mathcal{M}_{Y_{\mu}})$ 
be the minimal subfield of $\overline{\mathcal{M}_X}$ 
containing $i_{Y_{\lambda}}(\mathcal{M}_{Y_{\lambda}})$ and
$i_{Y_{\mu}}(\mathcal{M}_{Y_{\mu}})$. We note that
$\mathcal{L}$ is a finite Galois extension of
$\mathcal{M}_X$ by its definition. Let $g:Z\rightarrow X$ be the normalization of $X$ in $\mathcal{L}$ as explained in Remark~\ref{rk:normalization}. By the construction, we get an
object $(Z, g, i_Z)\in \mbox{Ob }\mathcal{I}_1^m(X;D)^p$ dominating both of $(Y_{\lambda},f_{\lambda},i_{Y_{\lambda}})$ and $(Y_{\mu},f_{\mu},i_{Y_{\mu}})$ in $\mathcal{I}_1^m(X;D)^p$. Let
$(\tilde X,\pi,i_{\tilde X})\in \mbox{Ob }\mathcal{I}_1^{\dag}G(X;D)^p$ be a final object of $\mathcal{I}_1^{\dag}G(X;D)^p$. From the equality $i_{\tilde X}=i_Y\circ \varpi_Y(i_Y,i_{\tilde
X})^{\ast}$, we see that there exists a canonical embedding: 
\begin{eqnarray}\label{eq:embed}
\Phi_{i_{\tilde X}}:\mathcal{I}_1^{\dag}G(X;D)^p\rightarrow \mathcal{I}_1^{\dag}G(\tilde X;0)^p, 
\end{eqnarray} 
depending on the choice of pointings $i_{\tilde X}$ for $(\tilde X,\pi)\in \mbox{Ob }\mathcal{I}_1^{\dag}G(X;D)$, such that $\Phi_{i_{\tilde X}}((Y,f,i_Y))=(Y,\varpi_Y(i_Y,i_{\tilde X}),
i_Y)\in \mbox{Ob }\mathcal{I}_1^{\dag}G(\tilde X;0)^p$ for $(Y,f,i_Y)\in \mbox{Ob }\mathcal{I}_1^{\dag}G(X;D)^p$.
Since $\mathcal{I}_1^{\dag}G(\tilde X;0)^p$ is cofilterd as explained in Remark~\ref{rk:prill-fun}, there exists $(W,h,i_W)\in \mbox{Ob }\mathcal{I}_1^{\dag}G(\tilde X;0)^p$ which dominates
both of $(Y_{\lambda},\varpi_{Y_{\lambda}}(i_{Y_{\lambda}},i_{\tilde X}), i_{Y_{\lambda}})$ and $(Y_{\mu},\varpi_{Y_{\mu}}(i_{Y_{\mu}},
i_{\tilde X}), i_{Y_{\mu}})$ in 
$\mathcal{I}_1^{\dag}G(\tilde X;0)^p$. By the construction of
$(Z,g,i_Z)\in \mbox{Ob }\mathcal{I}_1^m(X;D)^p$, $i_Z$ factors into $i_W\circ \tau^{\ast}$, where $\tau^{\ast}:\mathcal{M}_Z\rightarrow \mathcal{M}_W$ is an injective homomorphism.
Let
$\tau$ be the induced morphism $\tau:W\rightarrow Z$. Then we see that 
$h$ factors into $\varpi_Z(i_Z,i_{\tilde X})\circ \tau$. Since $h$ is \'etale in codimension one, hence so is $\varpi_Z(i_Z,i_{\tilde X})$. Thus we conclude that 
$(Z, g,i_Z)\in \mbox{Ob }\mathcal{I}_1^{\dag}G(X;D)^p$ and consequently, 
$\mathcal{I}_1^{\dag}G(X;D)^p$ is cofilterd. As for second statement, take any 
$(Y,f,i_Y)\in \mbox{Ob }\mathcal{I}_1^{\dag}(X;D)^p$ and let $\{i_Y^{(k)}|k=1,2,\dots n\}$ be all the pointings for 
$(Y,f)\in \mbox{Ob }\mathcal{I}_1^{\dag}(X;D)$. Let 
$\mathcal{L}$ be the minimal subfield of $\overline{\mathcal{M}_X}$ containing all the subfields $i_Y^{(1)}(\mathcal{M}_Y)$, \dots, $i_Y^{(1)}(\mathcal{M}_Y)$. Since $\mathcal{L}$ is a finite Galois
extension of $\mathcal{M}_X$ by its construction, we have an object $(Z,g,i_Z)\in \mbox{Ob }\mathcal{I}_1^mG(X;D)^p$ dominating all the $(Y,f,i_Y^{(1)})$, \dots, $(Y,f,i_Y^{(n)})\in \mbox{Ob
}\mathcal{I}_1^{\dag}(X;D)^p$, where $g:Z\rightarrow X$ is the normalization of $X$ in
$\mathcal{L}$. In the same way as in the previous argument, we conclude that
$(Z,g,i_Z)\in \mbox{Ob }\mathcal{I}_1^{\dag}G(X;D)^p$.
\hfill\bsquare
\begin{theorem}\label{th:comparison} For any $D\in \mbox{\em Div}_{\mbox{\boldmath $Q$}}X$, $\hat\pi_{1,X,p}^{{\rm loc}}[D]$ is isomorphic
to  $\overline\pi_{1,{\rm loc}}^{D^{\vee}}(\mbox{\em Reg }X)^{\wedge}$ in the category of profinite groups, where $D^{\vee}$ was defined as in Notation~\ref{notation:Dvee}.
\end{theorem}
{\it Proof.} By Proposition~\ref{pr:top-int} and Remark~\ref{rk:top-int}, we only have to show that $\mathcal{I}_1^{\dag}G(X;D)^p$ is cofinal in $\mbox{Int}^{\dag}G(X;D)^p$ 
(see \cite{sga4}, Expos\'e I, Proposition 8.1.3 or \cite{artin-mazur}, Appendix, Corollary (2.5)). Choose any object $(Y,f,i_Y)\in \mbox{Ob }\mbox{Int}^{\dag}G(X;D)^p$ and let
$\pi_Y:\tilde Y\rightarrow Y$ be the canonical cover with respect to $f^{\ast}D$. We can choose a pointing
$i_{\tilde Y}$ so that a triple $(\tilde Y,\tilde f,i_{\tilde Y})$ becomes an object in $\mbox{Int}^{\dag}(X;D)^p$ dominating $(Y,f,i_Y)$. From Remark~\ref{rk:lift}, we deduce that
$(\tilde Y,\tilde f,i_{\tilde Y})\in \mbox{Ob }\mathcal{I}_1^{\dag}G(X;D)^p$ by its construction.
\hfill\bsquare
\begin{rk}\label{rk:comparison}{\em In particular, assume that $D$ is integral, that is, $D\in \mbox{Div}_{\mbox{\boldmath $Q$}}X\cap\mbox{Weil }X$. Then Theorem~\ref{th:comparison}
says that 
$\hat\pi_{1,X,p}^{{\rm loc}}[D]$ is isomorphic to 
$\hat\pi_1^{{\rm loc}}(\mbox{Reg }X)$ as a profinite group.
}\end{rk}
\section{Universal index one covers}
Let $U$ be a Prill's good neighbourhood with regard to $\mbox{Sing }X$ and $\{U_{\lambda}\}_{\lambda\in \Lambda}$ its associated neighbourhood basis. Take any $(Y,f, i_Y)\in \mbox{Ob
}\mathcal{I}_1^{\dag}(X;0)^p$ and put
$(V^{-}, f^-, y):=\mathcal{P}^{-1}(Y,f, i_Y)\in FT(U^{-})$, where $\mathcal{P}$ is a Prill functor. Let $f:V\rightarrow U$ be the extended finite cover of $f^-$.
By Lemma~\ref{lm:lifting of homotopy}, $V$ is a Prill's good neighbourhood with regard to $f^{-1}(\mbox{Sing }X)$ with $\{V_{\lambda}\}_{\lambda\in \Lambda}$ being its associated
neighbourhood basis. Thus we have 
$\pi_1(V^-)=\mbox{projlim}_{q\in \mathcal{V};{\rm open}}\pi_1(\mathcal{V}\setminus f^{-1}(\mbox{Sing }X))=
\mbox{projlim}_{q\in \mathcal{V};{\rm open}}\pi_1(\mbox{Reg }\mathcal{V})=\pi_1^{{\rm loc}}(\mbox{Reg }Y)$ since $\mbox{Reg }\mathcal{V}\cap f^{-1}(\mbox{Sing }X)$ is a closed analytic subspace of codimension at least two in
$\mbox{Reg }\mathcal{V}$. (see, for example, 
\cite{prill}, III, Corollary 2.). In particular, we see that $(Y,f, i_Y)\in \mbox{Ob
}\mathcal{I}_1^{\dag}(X;0)^p$ is an initial object of $\mathcal{I}_1^{\dag}(X;0)^p$ if and only if $\hat\pi_1^{{\rm loc}}(\mbox{Reg }Y)=\{1\}$.
\begin{df}\label{df:u-index one cover}{\em A strict index one Galois cover
$\pi^{\dag}:X^{\dag}\rightarrow X$ with respect to $D$ is called an {\it algebraic 
universal index one cover} of a pair $(X,D)$, or abbreviated, a {\it universal index one cover} with respect to 
$D$ if $\hat\pi_1^{{\rm loc}}(\mbox{Reg }X^{\dag})=\{1\}$.
}\end{df}
\begin{rk}{\em Singularity with trivial local algebraic fundamental group is quite
restrictive one. For example, $\hat\pi_1^{{\rm loc}}(\mbox{Reg}X)=\{1\}$ implies $\mbox{Div}_{\mbox{\boldmath $Q$}}X\cap \mbox{Weil }X=\mbox{Div }X$. Moreoever, if we assume, in
addition, that $(X,p)$ is analytically $\mbox{\boldmath
$Q$}$-factorial, then $\mathcal{O}_X$ is factorial (see, for example, \cite{brieskorn}, Satz 1.4).}
\end{rk}
\begin{pr}\label{pr:fcc}  There exists the universal index one cover of $X$ with respect to $D$ if and only if $\hat\pi_1^{{\rm loc}}(\mbox{\em Reg }\tilde X)$ is finite.
\end{pr}
{\it Proof.} Assume that $\hat\pi_1^{{\rm loc}}(\mbox{Reg }\tilde X)$ is finite and take a final object $(\tilde X, \pi,i_{\tilde X})\in \mbox{Ob }\mathcal{I}_1^{\dag}G(X;D)^p$.
By the assumption, $\mathcal{I}_1^{\dag}(\tilde X;0)^p$ has an initial object
$(Y,f,i_Y)\in
\mbox{Ob }\mathcal{I}_1^{\dag}G(\tilde X;0)^p$ such that 
$\hat\pi_1^{{\rm loc}}(\mbox{Reg }Y)=\{1\}$. We note that $i_Y$ is also a pointing for $(Y,\pi\circ f)\in \mbox{Ob }\mathcal{I}_1^m(X;D)$ since we have $i_Y\circ
f^{\ast}\circ\pi^{\ast}=i_{\tilde X}\circ \pi^{\ast}=i_X$. Consider an object
$(Y,\pi\circ f,i_Y)\in
\mbox{Ob }\mathcal{I}_1^m(X;D)^p$. By the argument in Remark~\ref{rk:lift}, (1), we see that
$\pi\circ f$ is Galois. Since $\varpi_Y(i_Y,i_{\tilde X})=f$ is \'etale in codimension one, we conclude that $(Y,\pi\circ f,i_Y)\in \mbox{Ob }\mathcal{I}_1^{\dag}G(X;D)^p$. Conversely,
assume that there exists a pointed universal index one cover $(X^{\dag},\pi^{\dag},i_{X^{\dag}})\in \mbox{Ob }\mathcal{I}_1^{\dag}G(X;D)^p$ with respect to $D$. Then 
$\Phi_{i_{\tilde X}}((X^{\dag},\pi^{\dag},i_{X^{\dag}}))\in \mbox{Ob }\mathcal{I}_1^{\dag}G(\tilde X;0)^p$ is an initial object of $\mbox{Ob }\mathcal{I}_1^{\dag}(\tilde X;0)^p$, hence
$\hat\pi_1^{{\rm loc}}(\mbox{Reg }\tilde X)$ is finite.
\hfill\bsquare
\begin{df}[\cite{shokurov-cp}]{\em {\bf Q}-divisor $\Delta$ is called a 
{\it standard {\bf Q}-boundary} if $\mbox{mult}_{\Gamma}\Delta\in \{(n-1)/n|n\in \mbox{\boldmath $N$}\cup\{\infty\}\}$ for any prime divisor $\Gamma$. 
}\end{df}
\begin{rk}\label{rk:s-w}{\em Assume that $(X,\Delta)$ is purely log terminal, where $\Delta$ is a standard {\bf Q}-boundary. Then $(\tilde X,\Delta_{\tilde X})$ is known to be
canonical, hence $\tilde X$ has only canonical singularity if we assume that $\lfloor\Delta\rfloor=0$ or
$\tilde X$ is {\bf Q}-Gorenstein. Thus if $\mbox{dim } X\leq 3$, then $\hat\pi_1^{{\rm loc}}(\mbox{Reg }\tilde X)$ is finite by
\cite{shepherd-wilson}, Theorem 3.6. }
\end{rk}
\begin{pr}[c.f., \cite{sga1}, Expos\'e IX, Remark 5.8]\label{pr:fund.ext seq} There exists the following exact sequence in the category of profinite groups $:$
\begin{eqnarray}\label{eq:g-extseq}
\{1\}\longrightarrow\hat\pi_1^{{\rm loc}}(\mbox{\em Reg }\tilde X)\longrightarrow \hat\pi_{1,X,p}^{{\rm loc}}[D]\longrightarrow \mbox{\em Gal }(\tilde X/X)\simeq
\mbox{\boldmath $Z$}/r\mbox{\boldmath $Z$}\longrightarrow
\{1\}.
\end{eqnarray}
\end{pr}
{\it Proof. } Recall that we have a canonical embedding
$\Phi_{i_{\tilde X}}:\mathcal{I}_1^{\dag}(X;D)^p\rightarrow \mathcal{I}_1^{\dag}(\tilde X;0)^p$ as in (\ref{eq:embed}). Since we have the exact sequence:
$$
\{1\}\longrightarrow\mbox{projlim}_{(Y,f,i_Y)\in {\rm Ob }\mathcal{I}_1^{\dag}G(X;D)^p}\mbox{Gal }(Y/\tilde X)\longrightarrow \hat\pi_{1,X,p}^{{\rm loc}}[D]\longrightarrow
\mbox{Gal }(\tilde X/X)\longrightarrow
\{1\},
$$
we only have to show that $\mathcal{I}_1^{\dag}G(X;D)^p$ is cofinal in $\mathcal{I}_1^{\dag}G(\tilde X;0)^p$ via the functor $\Phi_{i_{\tilde X}}$. Choose any object $(Y,f,i_Y)\in \mbox{Ob
}\mathcal{I}_1^{\dag}G(\tilde X;0)^p$. Then we see that $(Y,\pi\circ f,i_Y)\in \mbox{Ob }\mbox{Int}^{\dag}(X;D)^p$ since $\pi^{-1}(\mbox{Reg }X\setminus\mbox{Supp }B)\subset \mbox{Reg
}\tilde X$ and $\pi\circ f$ is \'etale over $\mbox{Reg }X\setminus\mbox{Supp }B$. By Proposition~\ref{pr:top-int}, There exists an object $(Z,g,i_Z)\in \mbox{Ob
}\mbox{Int}^{\dag}G(X;D)^p$ dominating the object $(Y,\pi\circ f,i_Y)$.  Since $\mathcal{I}_1^{\dag}G(X;D)^p$ is cofinal in $\mbox{Int}^{\dag}G(X;D)^p$ (see the proof of
Theorem~\ref{th:comparison}), $(Z,g,i_Z)$ is dominated by some object in $\mathcal{I}_1^{\dag}G(X;D)^p$. Thus we get the assertion.
\hfill\bsquare
\begin{co}\label{co:cor to fund.ext seq} A pointed universal index one cover $(X^{\dag},\pi^{\dag},i_{X^{\dag}})\in \mbox{\em Ob }\mathcal{I}_1^{\dag}G(X;D)^p$
 is an initial object, or equivalently, a limit
of $\mathcal{I}_1^{\dag}(X;D)^p$ and vice versa. In particular, a universal index one cover with respect to 
$D$ is unique up to isomorphisms over $X$ if it exists. 
\end{co}
{\it Proof. }Let $(\tilde X,\pi,i_{\tilde X})\in \mbox{Ob }\mathcal{I}_1^{\dag}G(X;D)^p$ be a final object of $\mathcal{I}_1^{\dag}(X;D)^p$. Since
$(X^{\dag},\varpi_{X^{\dag}}(i_{X^{\dag}},i_{\tilde X}),i_{X^{\dag}})\in \mbox{Ob }\mathcal{I}_1^{\dag}(\tilde X;0)^p$ is an initial object of $\mathcal{I}_1^{\dag}(\tilde X;0)^p$,
$(X^{\dag},\pi^{\dag},i_{X^{\dag}})\in \mbox{Ob }\mathcal{I}_1^{\dag}G(X;D)^p$ is also an initial object of
$\mathcal{I}_1^{\dag}(X;D)^p$. On the contrary, assume that there exists an initial object $(X^{\dag\prime},\pi^{\dag\prime},i_{X^{\dag\prime}})$ of $\mathcal{I}_1^{\dag}(X;D)^p$. By
Proposition\ref{pr:cofilterd}, we see that $(X^{\dag\prime},\pi^{\dag\prime},i_{X^{\dag\prime}})\in \mbox{Ob }\mathcal{I}_1^{\dag}G(X;D)^p$. Since $\hat\pi_{1,X,p}^{{\rm loc}}[D]$ is
finite, $\hat\pi_1^{{\rm loc}}(\mbox{Reg }\tilde X)$ is also finite by Proposition~\ref{pr:fund.ext seq}. Therefore there exists a pointed universal index one cover
$(X^{\dag},\pi^{\dag},i_{X^{\dag}})\in \mbox{Ob }\mathcal{I}_1^{\dag}G(X;D)^p$ by Proposition~\ref{pr:fcc} which is also an initial object of $\mathcal{I}_1^{\dag}(X;D)^p$ and hence
isomorphic to $(X^{\dag\prime},\pi^{\dag\prime},i_{X^{\dag\prime}})$. Thus we conclude that $\hat\pi_1^{{\rm loc}}(\mbox{Reg }X^{\dag\prime})=\{1\}$.
\hfill\bsquare
\\[2ex]
The following Lemma is an algebraic generalization of a Brieskorn's lemma which we shall not use but will be useful in the other context.
\begin{lm}[c.f., \cite{brieskorn}, Lemma 2.6]\label{lm:brieskorn's Lemma}
Let $f:(X,p)\rightarrow (Y,q)$ be a finite morphism between germs of normal complex analytic spaces $(X,p)$ and $(Y,q)$. Then for any $D\in \mbox{\rm Div}_{\mbox{\boldmath $Q$}}Y$,
there exists a canonical homomorphism $f_{\ast}:\hat\pi_{1,X,p}^{{\rm loc}}[f^{\ast}D]\rightarrow \hat\pi_{1,Y,q}^{{\rm loc}}[D]$ which satisfies 
$|\hat\pi_{1,Y,q}^{{\rm loc}}[D]:\mbox{\rm Im }f_{\ast}|\leq \deg f$. In particular, if $\hat\pi_{1,X,p}^{{\rm loc}}[f^{\ast}D]$ is finite, so is 
$\hat\pi_{1,Y,q}^{{\rm loc}}[D]$.
\end{lm}
{\it Proof. } For a given pointing $i_Y:\mathcal{M}_Y\rightarrow \overline{\mathcal{M}_Y}$, we choose a pointing $i_X:\mathcal{M}_X\rightarrow \overline{\mathcal{M}_X}=\overline{\mathcal{M}_Y}$ such that
$i_X\circ f^{\ast}=i_Y$. Take any $(Z,\alpha,i_Z)\in \mbox{Ob }\mathcal{I}_1^{\dag}G(Y;D)^p$. Let $\nu:W\rightarrow Y$ be the normalization of $Y$ in 
$i_X(\mathcal{M}_X)\lor i_Z(\mathcal{M}_Z)$. We note that there exist morphisms $\beta:W\rightarrow X$ and $\gamma:W\rightarrow Z$ such that $\alpha\circ\gamma=f\circ\beta=\nu$. Since
$\beta^{\ast}f^{\ast}D=\gamma^{\ast}\alpha^{\ast}D\in \mbox{Div }W$ and $\beta:W\rightarrow X$ is Galois (see, for example, \cite{nagata}, Theorem 3.6.3), we have $(W, \beta, i_W)\in
\mbox{Ob }\mathcal{I}_1^{m}G(X;f^{\ast}D)^p$ for a suitable pointing
$i_W$.  Let $(\tilde X,\pi_X, i_{\tilde X})$ ( resp. $(\tilde Y,\pi_Y, i_{\tilde Y})$ ) be a final object of $\mathcal{I}_1^{m}G(X;f^{\ast}D)^p$ 
( resp. $\mathcal{I}_1^{m}G(Y;D)^p$ ). Since $(\tilde X, f\circ \pi_X, i_{\tilde X})\in \mbox{Ob }\mathcal{I}_1^{m}(Y;D)^p$, there exists a morphism $\omega_{\tilde X}(i_{\tilde
X},i_{\tilde Y}):(\tilde X, f\circ \pi_X, i_{\tilde X})\rightarrow (\tilde Y,\pi_Y, i_{\tilde Y})$ in $\mathcal{I}_1^{m}(Y;D)^p$. Let $\nu^{\natural}:Y^{\natural}\rightarrow Y$ be the
normalization of $Y$ in $i_Z(\mathcal{M}_Z)\cap i_{\tilde X}(\mathcal{M}_{\tilde X})$. Since the induced finite morphism $\delta:Z\rightarrow Y^{\natural}$ is Galois, $i_Z(\mathcal{M}_Z)$
and
$i_{\tilde X}(\mathcal{M}_{\tilde X})$ are linearly disjoint over 
$i_{Y^{\natural}}(\mathcal{M}_{Y^{\natural}})$, that is, 
$i_Z(\mathcal{M}_Z)\otimes _{i_{Y^{\natural}}(\mathcal{M}_{Y^{\natural}})}i_{\tilde X}(\mathcal{M}_{\tilde X})\simeq i_W(\mathcal{M}_W)$ (see, for example, \cite{nagata}, Exercise 4.2.3). Let
$\eta\in W$ be the generic point of a prime divisor on $W$ and $\xi\in Z$ (resp. $\tilde\xi\in \tilde X$, resp. $\xi^{\natural}\in Y^{\natural}$) be its image on $Z$ (resp. $\tilde
X$, resp. $Y^{\natural}$). Consider the canonical morphism 
$\kappa:i_Z(\mathcal{O}_{Z,\xi})\otimes _{i_{Y^{\natural}}(\mathcal{O}_{Y^{\natural},\xi^{\natural}})} i_{\tilde X}(\mathcal{O}_{\tilde X, \tilde \xi})\rightarrow i_W(\mathcal{O}_{W,\eta})$ and put
$S:=i_{Y^{\natural}}(\mathcal{O}_{Y^{\natural}})\setminus \{0\}$. We note that since $i_Z(\mathcal{O}_{Z,\xi})$ is flat over $i_{Y^{\natural}}(\mathcal{O}_{Y^{\natural},\xi^{\natural}})$ by our construction,  
$i_Z(\mathcal{O}_{Z,\xi})\otimes _{i_{Y^{\natural}}(\mathcal{O}_{Y^{\natural},\xi^{\natural}})} i_{\tilde X}(\mathcal{O}_{\tilde X, \tilde \xi})$ is a free $i_{\tilde X}(\mathcal{O}_{\tilde X, \tilde \xi})$-module, in particular, a torsion free
$i_{Y^{\natural}}(\mathcal{O}_{Y^{\natural},\xi^{\natural}})$-module. Since 
$S^{-1}\kappa:S^{-1}(i_Z(\mathcal{O}_{Z,\xi})\otimes _{i_{Y^{\natural}}(\mathcal{O}_{Y^{\natural},\xi^{\natural}})} 
i_{\tilde X}(\mathcal{O}_{\tilde X, \tilde \xi}))\simeq S^{-1}(i_Z(\mathcal{O}_{Z,\xi}))\otimes _{i_{Y^{\natural}}(\mathcal{M}_{Y^{\natural}})} S^{-1}(i_{\tilde X}(\mathcal{O}_{\tilde X, \tilde
\xi}))\rightarrow i_W(\mathcal{M}_W)$ is injective by the previous argument, so is $\kappa$, hence, in particular, $\mbox{Im }\kappa$ is a normal subring of $i_W(\mathcal{O}_{W,\eta})$
whose total quotient ring coincides $i_W(\mathcal{M}_W)$ which implies that $\mbox{Im }\kappa=i_W(\mathcal{M}_W)$. Thus we conclude that $\kappa$ is an isomorphism and that $i_W(\mathcal{O}_{W,\eta})$ is flat over $i_{\tilde X}(\mathcal{O}_{\tilde X,\tilde \xi})$, which implies that $(W,\beta,i_W)\in \mbox{Ob }\mathcal{I}_1^{\dag}G(X;f^{\ast}D)^p$. The canonical inclusion
$\mbox{Gal }(W/X)\simeq \mbox{Gal }(i_Z(\mathcal{M}_Z)/i_Z(\mathcal{M}_Z)\cap i_X(\mathcal{M}_X))\rightarrow \mbox{Gal }(Z/Y)$ induces a homomorphism $f_{\ast}:\hat\pi_{1,X,p}^{{\rm
loc}}[f^{\ast}D]\rightarrow \hat\pi_{1,Y,q}^{{\rm loc}}[D]$. Since we have $[\mbox{Gal }(Z/Y):\mbox{Gal }(W/X)]=[i_Z(\mathcal{M}_Z)\cap i_X(\mathcal{M}_X):i_Y(\mathcal{M}_Y)]\leq \deg f$, we
get the assertion (see also \cite{bourbaki}, \S 7.1, Corollaire 3).
\hfill\bsquare
\section{Lefshetz type theorem for $D$-local algebraic fundamental groups}
The aim of this section is to state and prove the Lefshetz type theorem for $D$-local algebraic fundamental groups.  
\begin{lm}[c.f., \cite{kawamata crep.}, Corollary 10.8]\label{lm:inv} Take any $D\in \mbox{\em Div}_{\mbox{\boldmath $Q$}}X\cap \mbox{\em Weil }X$ and let $\pi:\tilde X\rightarrow X$ be
the canonical cover with respect to D. Assume that there exists a normal prime divisor $\Gamma$ passing through
$p\in X$ such that the following three conditions hold. 
\begin{enumerate}
\item \ $\tilde \Gamma:=\pi^{-1}\Gamma$ is normal, 
\item \ $\Gamma$ does not contained in $\mbox{\em Supp }D$, 
\item \ there exists an analytic closed subset
$\Sigma \subset X$ with $\mbox{\em codim}_X\Sigma\geq 2$ and $\mbox{\em
codim}_{\Gamma}(\Sigma\cap \Gamma)\geq 2$ such that $D|_U$ is Cartier and
$D_{\Gamma}:=j_{\Gamma}^{\ast}D\in \mbox{\em Div }\Gamma$, where $U:=X\setminus \Sigma$
and $j_{\Gamma}:\Gamma_0:=\Gamma\setminus \Sigma\rightarrow \Gamma$ is the natural embedding.  
\end{enumerate} Then $D\in \mbox{\em Div }X$.
\end{lm}
{\it Proof.} Consider the exact sequence: 
$$
0\longrightarrow \mathcal{O}_{\tilde X}(\pi^{\ast}D-\tilde\Gamma)\longrightarrow \mathcal{O}_{\tilde X}(\pi^{\ast}D)\longrightarrow 
\mathcal{O}_{\tilde\Gamma}(\pi^{\ast}D)\longrightarrow 0. 
$$
Note that $\mathcal{O}_{\tilde \Gamma}(\pi^{\ast}D)$ and $\pi^{\ast}\mathcal{O}_{\Gamma}(D_{\Gamma})$ are both invertible and coincide on $\tilde\Gamma\setminus \pi^{-1}(\Sigma)$, hence we
have
$\mathcal{O}_{\tilde
\Gamma}(\pi^{\ast}D)=\pi^{\ast}\mathcal{O}_{\Gamma}(D_{\Gamma})$ by the normality of $\tilde\Gamma$. Since
$\pi_{\ast}^G:=\Gamma_X^G\circ\pi_{\ast}$ is an exact functor, where $G:=\mbox{Gal }(\tilde X/X)$, the above exact sequence induces a surjective map 
$$
\alpha:\mathcal{O}_X(D)=\pi^G_{\ast}\mathcal{O}_{\tilde X}(\pi^{\ast}D)\rightarrow \pi_{\ast}^G\mathcal{O}_{\tilde\Gamma}(\pi^{\ast}D)=
\mathcal{O}_{\Gamma}(D_{\Gamma})\otimes(\pi|_{\tilde\Gamma})_{\ast}^G\mathcal{O}_{\tilde\Gamma}=\mathcal{O}_{\Gamma}(D_{\Gamma}).
$$ Take $\varphi_{\Gamma}\in \mathcal{M}_{\Gamma}$ such that $\mbox{div }\varphi_{\Gamma}=-D|_{\Gamma}$. By the above argument, we have
$\varphi\in\mathcal{M}_X$ such that $\alpha(\varphi)=\varphi_{\Gamma}$. Since $(D+\mbox{div
}\varphi)|_{\Gamma}=D_{\Gamma}+\mbox{div }\varphi_{\Gamma}=0$ and $D+\mbox{div }\varphi$
is {\bf Q}-Cartier, we deduce that $\Gamma\cap\mbox{Supp }(D+\mbox{div
}\varphi)=\emptyset$, that is, $D+\mbox{div }\varphi=0$, hence $D\in \mbox{Div }X$.
\hfill\bsquare
\begin{lm}[c.f., \cite{reid c}, Lemma 1.12, \cite{seidenberg}]\label{lm:seidenberg} 
Let $X$ be a normal complex analytic space embedded in some domain in $\mbox{\boldmath
$C$}^n$. Consider the hypersurfaces
$H_{\tau}$ on $\mbox{\boldmath $C$}^n$ parametrized by $\tau\in \mbox{\boldmath
$P$}^n$ which is defined by a linear equation
$\tau_0+\sum_{i=1}^n\tau_iz_i=0$, where $z_1$, $\dots$, $z_n$ is a complete coordinate system of $\mbox{\boldmath
$C$}^n$. Then there exist a non-empty open subset $U\subset \mbox{\boldmath
$P$}^n$ and a countable union $Z$ of closed analytic subsets of $U$ such that for any $\tau\in U\setminus Z$, $H_{\tau}\cap X$ is a normal hypersurface on $X$. 
\end{lm}
{\it Proof.} Take an analytic open subset $U\subset \mbox{\boldmath
$P$}^n$ such that for any $\tau\in U$, $\bar H_{\tau}:=H_{\tau}\cap X$ is non empty and $H_{\tau}$
does not contain $X$. Since the base point free linear system $\{\bar H_{\tau}\}_{\tau\in U}$ on $X$
induces a base point free linear system on
$\mbox{Reg }X$, we have $\mbox{Sing }\bar H_{\tau}\subset \mbox{Sing } X$ and 
$\mbox{codim }_{\bar
 H_{\tau}}\mbox{Sing }\bar H_{\tau}\geq 2$ for any
$\tau\in U\setminus Z$, where $Z$ is a countable union of closed analytic subsets of $U$ by Bertini's theorem. Moreover, we may assume that for $k=1,\dots, d-2$, $\bar H_{\tau}$ does
not contain any maximal dimensional components of $(\mbox{Sing }X)\cap S_{k+1}(\mathcal{O}_X)$, where
$d:=\dim X$ and $S_{k}(\mathcal{O}_X)$ is a closed analytic set consisting of points at which the 
profoundity of $\mathcal{O}_X$ does not
exceed $k$. Since we have
$(\mbox{Sing }\bar H_{\tau})\cap S_{k}(\mathcal{O}_{\bar H_{\tau}})\subset (\mbox{Sing }X)\cap
S_{k+1}(\mathcal{O}_X)$, we see that $\dim (\mbox{Sing }\bar H_{\tau})\cap S_{k}(\mathcal{O}_{\bar
H_{\tau}})\leq k-2$ for any
$k$, hence $\bar H_{\tau}$ is normal for any $\tau\in U\setminus Z$ by \cite{fischer},  2.27, Theorem.
\hfill\bsquare
\begin{rk}{\em Let $X$ be a normal Stein space. For
{\bf Q}-divisor $\Delta$ on $X$, let
$\mbox{Mult}_X(\Delta)\subset
\mbox{\boldmath $Q$}$ denote the subset consisting of all the multiplicities of $\Delta$ at prime
divisors on $X$. We note that for general normal hyperplanes $\bar H_{\tau}$, we have 
$\mbox{Mult}_{\bar H_{\tau}}(\mbox{Diff}_{\bar H_{\tau}}(\Delta))\subset \mbox{Mult}_X(\Delta)$.
}\end{rk}
To state the Lefshetz type theorem, we need to fix some sort of general conditions. We shall
consider the following conditions assuming $\dim X\geq 2$. 
\begin{itemize}
\item[$(M1)$]\ $\Delta$ is a standard {\bf Q}-boundary. 
\item[$(M2)$]\ $(X,\Delta)$ is divisorially log terminal.
\item[$(M2)^{\ast}$]\ $(M2)^{\alpha}$ $(X,\Delta)$ is divisorially log terminal and $\{\Delta\}=0$ or $(M2)^{\beta}$ $(X,\Delta)$ is purely log terminal.
\item[$(M3)$]\ There exists an irreducible component
$\Gamma$ of $\lfloor\Delta\rfloor$ passing through $p\in X$ such that $K_X+\Gamma$ is {\bf Q}-Cartier.
\end{itemize}
\begin{rk} {\em $(M2)^{\ast}$ is a slightly stronger condition than $(M2)$.}
\end{rk}
\begin{pr}\label{pr:ind} Assume the conditions $(M1)$, $(M2)$ and $(M3)$. Then 
$$
\mbox{\em ind}_p(K_X+\Delta)=\mbox{\em ind}_p(K_{\Gamma}+\mbox{\em Diff}_{\Gamma}(\Delta-\Gamma)).
$$
\end{pr}
{\it Proof.} Put $r_{\Gamma}:=\mbox{ind}_p(K_{\Gamma}+\mbox{Diff}_{\Gamma}(\Delta-\Gamma))$. 
Firstly, we note that $(X,\Gamma)$ is purely log terminal and that $\Gamma\cap \mbox{Supp }(\Delta-\Gamma)$ is purely one codimensional in $\Gamma$ since $\Delta-\Gamma$ is 
{\bf
Q}-Cartier by the conditions $(M2)$ and $(M3)$. We show that
$r_{\Gamma}(K_X+\Delta)$ is an integral divisor on
$X$ and is Cartier at general points of any prime divisors on $\Gamma$. By taking general
hyperplane sections, we only have to check that if $\dim X=2$, then $r_{\Gamma}(K_X+\Delta)$ is
Cartier. This can be checked by the classification of log canonical singularities with a standard {\bf Q}-boundary due to S. Nakamura
(see, \S 3.1 or 
\cite{kobayashi}, Theorem 3.1), but we can also argue in this way as follows. We note that $p\in X$ is a cyclic quotient singular point with the order, say, 
$n$ by the condition $(M2)$. If $\Gamma\cap(\lfloor\Delta\rfloor-\Gamma)\neq \emptyset$, then $X$ is smooth, hence this case is
trivial. Assume that $\Gamma\cap(\lfloor\Delta\rfloor-\Gamma)=\emptyset$. Since we can see that $(X,\Delta)$ is purely log terminal in this
case from the condition $(M2)$, we can write $\Delta=\Gamma+d\Xi$ for a prime divisor $\Xi$ such that
$(\Gamma,\Xi)_p=1$ and for some
$d=(l-1)/l$, where $l$ is a natural number and we have $\mbox{mult}_p\mbox{Diff}_{\Gamma}(\Delta-\Gamma)=(nl-1)/(nl)$, as in \cite{shokurov-cp}, Lemma
2.25, which implies $r_{\Gamma}(K_X+\Delta)\in \mbox{Div }X$. Going back to the general case, we see that
$D:=K_X+\Delta\in \mbox{Div}_{\mbox{\boldmath $Q$}}X$ and $\Gamma$ satisfies the conditions in Lemma~\ref{lm:inv} using
\cite{shokurov}, Corollary 2.2 and Lemma 3.6, hence we conclude that $r_{\Gamma}(K_X+\Delta)\in \mbox{Div }X$. 
\hfill\bsquare
\begin{rk} {\em We note that $\mbox{Diff}_{\Gamma}(\Delta-\Gamma)$ is also a standard {\bf Q}-boundary, 
since $\mbox{Diff}_{\tilde\Gamma}((\Delta-\Gamma)_{\tilde X})$ is a {\bf Q}-boundary (see \cite{shokurov}, (2.4.1)).
}\end{rk}
\begin{ex}{\em Let $X$ be the germ of $\mbox{\boldmath $C$}^2$ at the
origin and put $\Gamma:=\mbox{div }z$ and  
$\Delta:=\mbox{div }z+(1/n)\mbox{div }w+(1/n)\mbox{div }(z+w)$, where $(z,w)$ is a system of coordinates and $n\in \mbox{\boldmath $N$}$. Then we have $\mbox{ind}_0(K_X+\Delta)=n$ while
$\mbox{ind}_0(K_{\Gamma}+\mbox{Diff}_{\Gamma}(\Delta-\Gamma))=n/2$ (resp. $n$) if
$n$ is even (resp. if $n$ is odd), which explains why we need the assumptions in Proposition~\ref{pr:ind}. }
\end{ex}
A directed set $(\Lambda, \geq)$ naturally forms a cofilterd projective system assuming that for $\lambda$, $\mu\in \Lambda$, 
$\mbox{Card Hom}_{\Lambda}(\lambda,\mu)=1$ if and only if $\lambda\geq \mu$. We call this projective system $\Lambda$ a {\it cofilterd index
projective system}. Let us recall the following basic result (see, for example, \cite{ribes}).
\begin{lm}\label{lm:topgrp} Let $\phi:\Lambda^{\prime}\rightarrow \Lambda$ be a covariant functor between cofilterd index
projective systems and $G:\Lambda\rightarrow \mbox{\em (Top. groups)}$, $H:\Lambda^{\prime}\rightarrow \mbox{\em (Top. groups)}$ be
two covariant functors to the category of topological groups. Assume that the following three conditions $(a)$, $(b)$ and $(c)$
hold.
\begin{itemize}
\item[$(a)$ ] $G_{\lambda}:=G(\lambda)$ and
$H_{\lambda^{\prime}}:=H(\lambda^{\prime})$ are compact for any $\lambda\in \mbox{\em Ob }\Lambda$ and
$\lambda^{\prime}\in\mbox{\em Ob }\Lambda^{\prime}$.
\item[$(b)$ ] $G(\lambda\rightarrow \mu)$ and $H(\lambda^{\prime}\rightarrow \mu^{\prime})$ are all surjective.
\item[$(c)$ ] There exists a natural transformation $\Psi:G\circ\phi\rightarrow H$ such that
$\Psi(\lambda^{\prime}):G_{\phi(\lambda^{\prime})}\rightarrow H_{\lambda^{\prime}}$ are surjective for any 
$\lambda^{\prime}\in\mbox{\em Ob
}\Lambda^{\prime}$.
\end{itemize}
Then there exists a canonical surjective morphism in \mbox{\em (Top. groups)} $:$
$$
\psi:\mbox{\em projlim }_{\lambda\in {\rm Ob }\Lambda}G_{\lambda}\rightarrow \mbox{\em projlim }_{\lambda^{\prime}\in {\rm Ob
}\Lambda^{\prime}}H_{\lambda^{\prime}}.
$$
\end{lm}
\vskip 5mm
Let $\Delta$ be a {\bf Q}-divisor on $X$ such that $K_X+\Delta$ is {\bf Q}-Cartier. In what follows, we put 
$$
\mathcal{I}^{\dag(m)}_1(G)(X,\Delta)^{(p)}:=\mathcal{I}^{\dag(m)}_1(G)(X;K_X+\Delta)^{(p)}
$$ and 
$$
\hat\pi_{1,X,p}^{{\rm loc}}<\Delta>:=\hat\pi_{1,X,p}^{{\rm loc}}[K_X+\Delta].
$$
\begin{theorem}\label{th:pi1surj} Assume the conditions $(M1)$, $(M2)^{\ast}$ and $(M3)$. Then there exists a canonical continuous 
surjective homomorphism
$$
\psi_{\Gamma}:\hat\pi^{{\rm loc}}_{1,\Gamma,p}<\mbox{\em
Diff}_{\Gamma}(\Delta-\Gamma)>\rightarrow \hat\pi^{{\rm loc}}_{1,X,p}<\Delta>. 
$$\end{theorem}
{\it Proof. }For 
any
$(Y,f)\in \mbox{Ob }\mathcal{I}_1^m(X,\Delta)$, $\Gamma$ and $\Gamma_Y$ are normal by
\cite{shokurov}, Lemma 3.6 and Corollary 2.2, hence they are irreducible since $f^{-1}(p)$ consists
of just one point. A canonical inclusion 
$\mathcal{O}_{\Gamma}\rightarrow  \mbox{injlim}_{(Y,f,i_Y)\in {\rm Ob }\mathcal{I}_1^m(X,\Delta)^p}\mathcal{O}_{\Gamma_Y}$ extends to an inclusion 
$i_{\Gamma}:\mathcal{O}_{\Gamma}\rightarrow \overline{\mathcal{M}_{\Gamma}}$ and we fix this $i_{\Gamma}$. Then we have a canonical functor
$
\phi^{(p)}_{\Gamma}:\mathcal{I}_1^m(X,\Delta)^{(p)}\rightarrow \mathcal{I}_1^m(\Gamma,\mbox{Diff}_{\Gamma}(\Delta-\Gamma))^{(p)}
$
such that $\phi((Y,f))=(\Gamma_Y,f_{\Gamma})$, where $f_{\Gamma}:=f|_{\Gamma_Y}$. Take any $(Y,f)\in \mbox{Ob }\mathcal{I}_1^{\dag}(X,\Delta)$. 
We shall show that $(\Gamma_Y,f_{\Gamma})\in \mbox{Ob }\mathcal{I}_1^{\dag}(\Gamma,\mbox{Diff}_{\Gamma}(\Delta-\Gamma))$.
Let $\tilde X\rightarrow X$ be the canonical cover with respect to $K_X+\Delta$. Note that $\Gamma_{\tilde X}$ is also normal and 
$\pi_{\Gamma}:=\pi|{\Gamma_Y}:\tilde\Gamma:=\Gamma_{\tilde X}\rightarrow \Gamma$ is the canonical cover with respect to $K_{\Gamma}+\mbox{Diff}_{\Gamma}(\Delta-\Gamma)$ by
Proposition~\ref{pr:ind}. Take any pointings $i_{\tilde \Gamma}$ and $i_{\Gamma_Y}$ for 
$(\tilde \Gamma,\pi_{\Gamma})$ and $(\Gamma_Y,f_{\Gamma})\in \mbox{Ob
}\mathcal{I}_1^m(\Gamma,\mbox{Diff}_{\Gamma}(\Delta-\Gamma))$. We note also that $\varpi_{\Gamma_Y}(i_{\Gamma_Y},i_{\tilde \Gamma})=\varpi_Y(i_Y,i_{\tilde X})|_{\Gamma_Y}$ for some
pointings $i_Y$ and $i_{\tilde X}$ of $(Y,f)$ and $(\tilde X,\pi)\in \mbox{Ob }\mathcal{I}_1^{\dag}(X,\Delta)$. $(\Gamma_Y,f_{\Gamma})\in \mbox{Ob }\mathcal{I}_1^{\dag}(\Gamma,\mbox{Diff}_{\Gamma}(\Delta-\Gamma))$. By the covering theorem
in \cite{endre}, $(\tilde X,\tilde \Delta)$ is divisorially log terminal of index one, which
implies that $\tilde X$ is smooth in codimension two, hence, in particular, we have $\mbox{codim}_{\tilde\Gamma}(\mbox{Sing }\tilde X\cap\tilde\Gamma)\geq 2$. 
Since $\varpi_Y(i_Y,i_{\tilde X})$ is \'etale over
$\mbox{Reg }\tilde X$, we conclude that $\varpi_{\Gamma_Y}(i_{\Gamma_Y},i_{\tilde \Gamma})$ is \'etale in codimension one and $(\Gamma_Y,f_{\Gamma})\in \mbox{Ob }\mathcal{I}_1^{\dag}(\Gamma,\mbox{Diff}_{\Gamma}(\Delta-\Gamma))$. In other words, $\phi^{(p)}$ induces a functor
$
\phi^{(p)}_{\Gamma}:\mathcal{I}_1^{\dag}(G)(X,\Delta)^{(p)}\rightarrow \mathcal{I}_1^{\dag}(G)(\Gamma,\mbox{Diff}_{\Gamma}(\Delta-\Gamma))^{(p)},
$
where we used the same notation $\phi^{(p)}$. Consider the two functors
$$
G_X:\mathcal{I}_1^{\dag}G(X,\Delta)^p\rightarrow \mbox{(Top. groups)}\mbox{ and } G_{\Gamma}:\mathcal{I}_1^{\dag}G(\Gamma,\mbox{Diff}_{\Gamma}(\Delta-\Gamma))^p\rightarrow \mbox{(Top.
groups)}
$$
such that $G_X((Y,f,i_Y))=\mbox{Gal }(Y/X)$ and $G_{\Gamma}((\Gamma^{\prime},g,i_{\Gamma^{\prime}}))=\mbox{Gal }(\Gamma^{\prime}/\Gamma)$. Since $f$ is \'etale over a general points of
$\Gamma$ for any $(Y,f)\in \mbox{Ob }\mathcal{I}_1^{\dag}G(X,\Delta)$, there exists a natural equivalence $\Psi_{\Gamma}:G_{\Gamma}\circ\phi^p\rightarrow G_X$, which induces the desired
surjection 
$\psi_{\Gamma}:\hat\pi^{{\rm loc}}_{1,\Gamma,p}<\mbox{Diff}_{\Gamma}(\Delta-\Gamma)>\rightarrow \hat\pi^{{\rm loc}}_{1,X,p}<\Delta>$ by Lemma~\ref{lm:topgrp}.
\hfill\bsquare
\begin{rk}{\em Assume the conditions $(M1)$, $(M2)^{\alpha}$ and $(M3)$. Then, combined with Remark\ref{rk:comparison}, Theorem~\ref{th:pi1surj} says that there exists a surjection 
$
\psi_{\Gamma}:\hat\pi^{{\rm loc}}_{1,\Gamma,p}<\mbox{Diff}_{\Gamma}(\Delta-\Gamma)>\rightarrow \hat\pi_1^{{\rm loc}}(\mbox{Reg }X).
$
For example, if $\dim X=4$, $\hat\pi^{{\rm loc}}_{1,\Gamma,p}<\mbox{Diff}_{\Gamma}(\Delta-\Gamma)>$ is finite under the assumptions as explained in Remark~\ref{rk:s-w}, hence so
is 
$\hat\pi_1^{{\rm loc}}(\mbox{Reg }X)$. 
}\end{rk}
\begin{co}\label{co:induction} Let notation and assumptions be as in Theorem~\ref{th:pi1surj}. Assume that the universal index one cover of $\Gamma$ with respect to
$K_{\Gamma}+\mbox{\em Diff}_{\Gamma}(\Delta-\Gamma)$ exists. Then there exists the universal index one cover of $X$ with respect to $K_X+\Delta$. Moreover, there exists the following
exact sequence $:$
\begin{eqnarray}\label{eq:canexseq}
\{1\}\longrightarrow \hat\pi_1^{{\rm loc}}(\mbox{\em Reg
}\Gamma_{X^{\dag}})\longrightarrow \hat\pi^{{\rm loc}}_{1,\Gamma,p}<\mbox{\em Diff}_{\Gamma}(\Delta-\Gamma)>\longrightarrow \hat\pi^{{\rm loc}}_{1,X,p}<\Delta>\longrightarrow
\{1\},
\end{eqnarray}
where $\pi^{\dag}:X^{\dag}\rightarrow X$ is the universal index one cover of $X$ with respect to $K_X+\Delta$.
\end{co}
{\it Proof.} The first assertion follows from Proposition~\ref{pr:fcc} and Proposition~\ref{pr:fund.ext seq}. As for the last statement, let
$\pi^{\dag}_{\Gamma}:\Gamma^{\dag}\rightarrow \Gamma$ be  the universal index one cover of $\Gamma$ with respect to
$K_{\Gamma}+\mbox{Diff}_{\Gamma}(\Delta-\Gamma)$. Then the induced morphism $\tau^{\dag}_{\Gamma}:\Gamma^{\dag}\rightarrow \Gamma_{X^{\dag}}$ is the universal index
one cover of $\Gamma_{X^{\dag}}$ since $\tau^{\dag}_{\Gamma}$ is \'etale in codimension one, which implies that $\mbox{Gal }(\Gamma^{\dag}/\Gamma_{X^{\dag}})\simeq
\hat\pi_1^{{\rm loc}}(\mbox{Reg }\Gamma_{X^{\dag}})$, hence we obtain the desired exact sequence. 
\hfill\bsquare
\begin{rk}\label{rk:3-plt}{\em Let notation be as above. Assume that $(X,p)$ is a three dimensional {\bf Q}-Gorenstein singularity and that $(X,\Gamma)$ is purely log terminal with
$\mbox{Sing X}\subset
\Gamma$. Then we see that
$(\tilde X,\tilde p)$ has only terminal singularities and that  
$(X^{\dag},p^{\dag})$ is an isolated compound Du Val singularity (see, \cite{milnor}, Theorem 5.2). We also note that $\Gamma^{\dag}$ is smooth and that $\Gamma_{X^{\dag}}\in
|-K_{X^{\dag}}|$ is a Du Val element.  Moreover, the above
exact sequence (\ref{eq:canexseq}) reduces to the following exact sequence:
\begin{eqnarray}\label{eq:extseq3}
\{1\}\longrightarrow \hat\pi_1^{{\rm loc}}(\mbox{Reg
}\Gamma_{X^{\dag}})\longrightarrow \hat\pi^{{\rm loc}}_{1,\Gamma,p}<\mbox{Diff}_{\Gamma}(0)>\longrightarrow \hat\pi^{{\rm loc}}_1(\mbox{Reg }X)\longrightarrow \{1\},
\end{eqnarray}
which enables us to calculate the algebraic local fundamental group of $(X,p)$, since $\hat\pi_1^{{\rm loc}}(\mbox{Reg
}\Gamma_{X^{\dag}})$ and $\hat\pi^{{\rm loc}}_{1,\Gamma,p}<\mbox{Diff}_{\Gamma}(0)>$ have faithful representations to the special unitary group $SU(2,\mbox{\boldmath $C$})$ and the
unitary group
$U(2,\mbox{\boldmath
$C$})$ respectively, both of which are classified. It is important to determine the pair $(X^{\dag}, \hat\pi^{{\rm loc}}_1(\mbox{Reg }X))$ which will lead us to the classification
$3$-dimensional purely log terminal singularities. }\end{rk}
\begin{ex}{\em Let notation and assumptions be as in Remark~\ref{rk:3-plt}. Assume that $(\Gamma,p)\simeq (\mbox{\boldmath $C$}^2,0)$ and
$\mbox{Diff}_{\Gamma}(0)=(1/2)\mbox{div}(z^2+w^n)$ ($n\geq 2$), where $(z,w)$ is a system of coordinate of $\Gamma$. Moreover, assume that $(\tilde X,\tilde p)$ has only quotient
terminal singularity. Then we can deduce that $n=2$ by using our theory in such a way as follows. It can be easily checked that $(X^{\dag},p^{\dag})$ and
$(\Gamma^{\dag},p^{\dag})$ are both smooth and that $\hat\pi^{{\rm loc}}_{1,\Gamma,p}<\mbox{Diff}_{\Gamma}(0)>\simeq \hat\pi^{{\rm loc}}_1(\mbox{Reg }X)\simeq G$, where $G$ is the
dihedral group of the order $2n$ (see also \cite{nakano}). Let $G=<a,b;a^n=1,b^2=1,b^{-1}ab=a^{-1}>$ be a presentation of $G$ and $\rho_{\Gamma}:G\rightarrow U(2, \mbox{\boldmath $C$})$
be a corresponding representation with respect to $(\Gamma,\mbox{Diff}_{\Gamma}(0))$ defined as follows. 
$$
\rho_{\Gamma}(a)=\left(
\begin{array}{cc}
e^{2\pi i/n} & 0\\
0 & e^{-2\pi i/n}
\end{array}
\right),
\quad
\rho_{\Gamma}(b)=
\left(
\begin{array}{cc}
0 & 1\\
1 & 0
\end{array}
\right).
$$
Let $\rho_X:G\rightarrow U(3,\mbox{\boldmath $C$})$ be a corresponding faithful representation with respect to $(X,p)$. Since $\Gamma^{\dag}\subset X^{\dag}$ is invariant under the
action of
$G$ through $\rho_X$, $\rho_X$ is equivalent to $\rho_{\Gamma}\oplus \chi$ for some character $\chi:G\rightarrow \mbox{\boldmath $C$}^{\times}$. Let $K$ be the kernel of the character 
$G\rightarrow \mbox{Aut }\mathcal{O}_{X^{\dag}}(K_{X^{\dag}}+\Gamma^{\dag})/m_{p^{\dag}}\mathcal{O}_{X^{\dag}}(K_{X^{\dag}}+\Gamma^{\dag})$ induced by $\rho_X$. Then we see that $K=\mbox{Ker
 det }\rho_{\Gamma}$. Since we have $\mbox{det }\rho_{\Gamma}(a)=1$ and $\mbox{det }\rho_{\Gamma}(b)=-1$, we have $<a>\subset K$ and $b\notin K$. Noting that we have  
$2=[G:<a>]=[G:K][K:<a>]$, we obtain $K=<a>$. We also note that we have $\mbox{ord }\chi(a)=n$ since $\tilde p\in\tilde X\simeq \mbox{\boldmath $C$}^3/K$ is isolated. On the other hand,
since we have
$\chi(b^{-1}ab)=\chi(a^{-1})$, we get $\chi(a)^2=1$. Thus we conclude that $n=2$. This result will be used in \cite{ohno}.
}\end{ex}
The results in this paper, especially Theorem~\ref{th:pi1surj}, are heavily dependent on the results such as Hironaka's resolution theorem and Kawamata-Viehweg's
vanishing theorem and so on which are valid only in the case of characteristic zero.  
\begin{prob}{\em  Construct our theory in the positive characteristic case.
}\end{prob}


\begin{thebibliography}{1}
\bibitem{abhyankar} S. S. Abhyankar, {\em Local Analytic Geometry}, Academic Press, New York-London, 1964.
\bibitem{artin-mazur} M. Artin and B. Mazur, {\em Etale Homotopy}, Lecture Notes in Math. {\bf 100}, Springer Heidelberg, 1969.
\bibitem{bourbaki} N. Bourbaki, {\em Topologie G\'en\'erale}, Hermann, Paris, 1971.
\bibitem{brieskorn}  E. Brieskorn, {\em Rational Singularit\"aten komplexer Fl\"achen}, Invent. Math. {\bf 4}, (1968), pp.336-358.
\bibitem{catanese} F. Catanese, {\em Automorphisms of Rational Double Points and Moduli Spaces of
Surfaces of General Type}, Compos. Math. {\bf 61}, (1987), pp. 81-102.
\bibitem{fischer} G. Fisher, {\em Complex Analytic Geometry}, Lecture Notes in Math. {\bf 538},
Springer-Verlag, Berlin-Heidelberg-New York, 1976.
\bibitem{grauert-remmert} H. Grauert and R. Remmert, {\em Komplexe R\"aume}, Math. Ann. {\bf 136}, (1958), pp. 245-318.
\bibitem{grothendieck} A. Grothendieck, {\em Techniques de Construction en G\'eom\'etrie Analytique}, S\'eminaire Henri Cartan,
13ieme ann\'ee (1960/1).
\bibitem{grothendieck-murre} A. Grothendieck and J.P. Murre, {\em The Tame Fundamental Group of a Formal Neighbourhood of a Divisor with Normal Crossings on a Scheme}, 
Lecture Notes in Math. {\bf 208},  Springer-Verlag, Berlin-Heidelberg-New York, 1970.
\bibitem{sga1} A. Grothendieck et. al., {\em Rev\^etements Etales et Groupe Fondamental},
Lecture Notes in Math. {\bf 224},  Springer-Verlag, Berlin-Heidelberg-New York, 1971. 
\bibitem{sga4} A. Grothendieck et. al., {\em T\'eorie des Topos et Cohomologie Etale des Sch\'emas},
Lecture Notes in Math. {\bf 269},  Springer-Verlag, Berlin-Heidelberg-New York, 1972. 
\bibitem{gunning-rossi} R. C. Gunning and H. Rossi, {\em Analytic Functions of Several Complex Variables}, Prentice-Hall, Inc., Englewood Cliffs, N. J. 1965.
\bibitem{hironaka} H. Hironaka, {\em Resolution of Singularities of an Algebraic Variety over a
Field of Characteristic Zero I, II}, Ann. of Math. {\bf 79}, (1964), pp. 109-326.
\bibitem{kato} M. Kato, {\em On the uniformizations of orbifolds}, Adv. Stud. in Pure Math. {\bf 9}, (1986), pp.149-172.
\bibitem{kawamata crep.} Y. Kawamata, 
{\em Crepant blowing-up of $3$-dimensional canonical singularities 
and its application to degeneration of surfaces,} Jour. of Ann. of Math. {\bf 127},
(1988), pp.93-163 
\bibitem{kns} R. Kobayashi, S. Nakamura and F. Sakai, {\em A Numerical Characterization of Ball
Quotients for Normal Surfaces with Branch Loci}, Proc. Japan Acad., {\bf 65}, Ser. A (1989),
pp.238-241.
\bibitem{kobayashi} R. Kobayashi, {\em Uniformization of Complex Surfaces,} 
Adv. Stud. in Pure Math. {\bf 18-II}, (1990), pp. 313-394.
\bibitem{utah} J. Koll\'ar et. al., {\em Flips  and
abundunce for algebraic threefolds,}
  Ast\'erisque {\bf 211}, A summer seminar at the university of Utah Salt Lake City. 1992.
\bibitem{matsumura} H. Matsumura, {\em Commutative Algebra}, W. A. Benjamin Co, New York, 1970.
\bibitem{milnor} J. Milnor, {\em Singular points of complex hypersurfaces}, Princeton Univ. Press, Princeton, New Jersey and Univ. of Tokyo Press, Tokyo, 1968.
\bibitem{mitchell} B. Mitchell, {\em Theory of Categories}, Academic Press, New York and London, 1965.
\bibitem{murre} J. Murre, {\em Lectures on an Introduction to Grothendiek's Theory of the Fundamental Group},
 Lecture Notes, Tata Institute of Fundamental Research, Bombay, 1967.
\bibitem{nagata} M. Nagata, {\em Field Theory}, Pure and Applied Mathematics, A series of Monographs and Textbooks,
Marcel Dekker, Inc.,  New York and Basel, 1977.
\bibitem{nakano} T. Nakano and K. Tamai, {\em On Some Maximal Galois Coverings over Affine and Projective Planes}, Osaka J. Math. {\bf 33} (1996), 347-364.
\bibitem{namba} M. Namba, {\em Branched coverings and algebraic functions}, Pitman Research Notes in Math. Series {\bf 161}, 
Longman Scientific and Technical, Harlow, John Wiley and Sons, New York, 1987.
\bibitem{ohno} K. Ohno, {\em The Euler Characteristic Formula for Logarithmic minimal Degenerations of surfaces with Kodaira Dimension Zero and its application to Calabi-Yau Threefolds with a pencil.}, Available for download at math.AG /0710.3641.
\bibitem{prill} D. Prill, {\em Local Classification of Quotients of Complex Manifolds by
Discontinuous Groups}, Duke Math. J. {\bf 34}, (1967) pp. 375-386.
\bibitem{reid c} M. Reid, {\em Canonical threefold}, G\'eom\'etry Alg\'ebrique Angers, A. Beauville
ed., Sijthoff \& Noordhoff, 1980, pp. 273-310.
\bibitem{ribes} L. Ribes, {\em Introduction of Profinite Groups and Galois Cohomology}, Queen's Papers in Pure and Applied
Mathematics-No. {\bf 24}, Queen's University, Kingston, Ontario, 1970.
\bibitem{seidenberg} A. Seidenberg, {\em The hyperplane sections of normal varieties}, Trans. of the
Amer. Math. Soc., {\bf 69}, No. 2, (1950), pp. 357-386.
\bibitem{seidenberg2} A. Seidenberg, {\em Saturation of an analytic ring}, Amer. Jour. of Math., {\bf 94}, (1972), pp. 424-430.
\bibitem{serre} J.P. Serre, {\em Topics in Galois Theory}, Jones and Bartlett Publishers, Boston-London, 1992.
\bibitem{shepherd-wilson} N. I. Shepherd-Barron and P. M. H. Wilson, {\em Singular $3$-folds with
Numerically Trivial First and Second Chern Classes}, J. Algebraic Geometry {\bf 3}, (1994),
pp.265-281.
\bibitem{shokurov} VV. Shokurov, {\em $3$-fold Log Flips,} Russian Acad. Sci. Izv.
Math. {\bf 40}, (1993), pp.95-202.
\bibitem{shokurov-cp} VV. Shokurov, {\em Complements on surfaces}, preprint, 1997.
\bibitem{steenrod} N. Steenrod, {\em The Topology of Fibre Bundles}, Princeton Math. Series {\bf 14}, Princeton UP, 1951.
\bibitem{endre} E. Szab\'o, {\em Divisorial log terminal singularities}, J. Math. Sci. Univ.
Tokyo {\bf 1}, (1994), pp. 631-639.
\bibitem{gabi} G. Teodosiu, {\em A Class of Analytic Coverings Ramified over $u^3=v^2$}, J. London
Math. Soc. (2) {\bf 38}, (1988), pp. 231-242.
\end{thebibliography}
\end{document}